\documentclass[aip, chaos, reprint]{revtex4-1}

 \usepackage[parfill]{parskip}
 \usepackage[utf8]{inputenc}
\usepackage{graphicx}
   \usepackage{epstopdf}
\usepackage{subfigure}

\usepackage{amsthm}

\usepackage{mathptmx}

 \usepackage{color}
 
 \newtheorem{theorem}{Theorem}
\newtheorem{lemma}[theorem]{Lemma}

\newcommand{\defeq}{:=}
\newcommand{\RR}{\mathbf{R}}
\newcommand{\ZZ}{\mathbf{Z}}
\newcommand{\NN}{\mathbf{N}}

\newcommand{\itin}{\mathcal{I}}

\newcommand{\textfrac}[2]{\textstyle \frac{#1}{#2} }

\usepackage{hyperref}

\renewcommand{\text}[1]{\mathrm{#1}}

\newcommand{\tM}{\tilde{M}}

\newcommand{\eqref}[1]{(\ref{#1})}   % for IOP journals

\graphicspath{ {figs/}} 
\newcommand{\limitset}{\Gamma_\lambda}

\newcommand{\elcantor}{C(\lambda)}
\newcommand{\angmap}{\varphi}
\newcommand{\ang}{\theta}
\newcommand{\elshift}{\Sigma_2^+}
\newcommand{\DiscontinuityRegion}{\text{Sing}(\lambda)}

\newcommand{\hinterval}[2]{J(#1,#2)}

\newcommand{\partition}{\mathcal{S}}
\newcommand{\MarkovianRegion}{\mathcal{J}_{0}}

\newcommand{\maximalinvariantset}{\Gamma_{\lambda}}

\begin{document}
%  \maketitle

\title[Strange attractors in non-elastic billiards]{Structure and evolution of strange attractors in non-elastic triangular billiards}
% \\
% \emph{Dedicated to the memory of Jorge Alejandro Hernández Tahuilán}}

%\author{Aubin Arroyo$^1$, Roberto Markarian$^2$ and David P.~Sanders$^3$}
\author{Aubin Arroyo}
\affiliation{ Instituto de Matem\'aticas, Unidad Cuernavaca, Universidad Nacional Aut\'onoma de M\'exico,
Apartado postal 273-3, Admon.\ 3, Cuernavaca, 62251 Morelos, Mexico}
\email{aubin@matcuer.unam.mx}

\author{Roberto Markarian}
\affiliation{Instituto de Matem\'atica y Estad\'istica (IMERL), Facultad de Ingenier\'ia, Universidad de la Rep\'ublica,
Montevideo 11300, Uruguay}
\email{roma@fing.edu.uy}
\author{David P.~Sanders}
\affiliation{Departamento de F\'isica, Facultad de Ciencias, Universidad Nacional Aut\'onoma de M\'exico, Ciudad Universitaria, 04510 M\'exico D.F., Mexico}
\email{dpsanders@ciencias.unam.mx}

\date{\today}
%
%\address{$^1$  Instituto de Matem\'aticas, Unidad Cuernavaca, Universidad Nacional Aut\'onoma de M\'exico,
%Apartado postal 273-3, Admon.\ 3, Cuernavaca, 62251 Morelos, Mexico}
%% \ead{aubin@matcuer.unam.mx}
%
%\address{$^2$ Instituto de Matem\'atica y Estad\'istica (IMERL), Facultad de Ingenier\'ia, Universidad de la Rep\'ublica,
%Montevideo 11300, Uruguay}
%% \ead{roma@fing.edu.uy}
%
%\address{$^3$ Departamento de F\'isica, Facultad de Ciencias, Universidad Nacional Aut\'onoma de M\'exico, Ciudad Universitaria, 04510 M\'exico D.F., Mexico}
%% \ead{dps@fciencias.unam.mx}

%\eads{\mailto{aubin@matcuer.unam.mx}, \mailto{roma@fing.edu.uy} and \mailto{dpsanders@ciencias.unam.mx} }

\date{\today}

\begin{abstract}
We study pinball billiard dynamics in an equilateral triangular table.
%, focusing on the equilateral triangle.
In such dynamics, collisions with the walls are non-elastic: the outgoing angle with the normal vector to the boundary is a uniform factor $\lambda < 1$ smaller than the incoming angle.
This leads to contraction in phase space for the discrete-time dynamics between consecutive collisions, and hence to attractors of zero Lebesgue measure, which are almost always
fractal strange attractors with chaotic dynamics, due to the presence of an expansion mechanism. We study the structure of these strange attractors and their evolution as the 
contraction parameter $\lambda$ is varied. For $\lambda \in (0, \frac{1}{3})$, we prove rigorously that the attractor has the structure of a Cantor set times an interval,
whereas for larger values of $\lambda$ the billiard dynamics gives rise to nonaccessible regions in phase space.
For $\lambda$ close to $1$, the attractor splits into three transitive components, the basins of attraction of which have fractal basin boundaries.
% ; this is related ??? to intermingled
% basins in symmetric systems.
\end{abstract}

\maketitle

% \clearpage
%  \tableofcontents

\begin{quotation}
Billiard models form an important class of dynamical systems, in which a particle collides with fixed boundary elements. 
The collision rule plays a crucial role in the dynamics:
for classical elastic collisions, phase-space volume is preserved, and the dynamics
is typical of Hamiltonian systems \cite{CM06}.
The recently-introduced \emph{pinball billiards}\cite{MPS10}, or \emph{pinbilliards}, 
instead have a non-elastic collision rule, reminiscent of the effect of the extra impulse in the normal direction received by a ball in a pinball machine when it hits an active bumper. 
The present authors previously studied\cite{AMS09} a rule in which the outgoing angle at a collision, with respect to the normal vector to the boundary, is shrunk by a factor $\lambda < 1$ compared to the incoming angle.  
In  convex billiard tables with focusing and flat boundary components, attracting periodic orbits, strange attractors, and bifurcation phenomena are then observed, depending on the boundary geometry.% ROBERTO QUITÓ: \cite{AMS09}.
When the curvature of a boundary element is non-zero, some hyperbolic (chaotic) properties are obtained, both in the elastic \cite{CM06} and non-elastic \cite{AMS09, MPS10} cases.

In order to make further progress in understanding the origin of the complicated behaviour induced by the non-elastic collision rule, in 
 this paper we restrict attention to pinball billiards in polygonal tables, specifically the equilateral triangle, which turns out to capture many of the main features of such systems. 
Here we again find hyperbolic (chaotic) strange attractors, whose structure evolves in complicated ways as $\lambda$ varies between $0$ and $1$. Due to the simplicity of the geometry, however, many of these structural changes may now fact be understood analytically.

\end{quotation}
%\section{Introduction}
%In particular, we will study a collision rule for which there is contraction of phase space, so that the dynamics exhibits phenomena typical 
%of dissipative systems. %, in particular periodic and chaotic (hyperbolic) attractors \cite{AMS09}.
%
%influenced by the presence both of the curved sections and of the discontinuities between different
%boundary components.

Polygonal billiard tables, in the classical case of elastic collisions,
exhibit dynamical properties which depend strongly on the angles of the table:
% The billiard map has a natural invariant measure with zero Lyapunov
% exponents. 
if the angles are rational multiples of $\pi$, then the dynamics are completely regular, and % for each 
only a few angles can be realised during a given trajectory. 
On the other hand, arbitrarily close to a given polygon is another
 (with vertices as close as desired), whose billiard map is ergodic \cite{KMS86}.
Furthermore, triangular billiards whose angles are all irrational with $\pi$ are conjectured to be mixing \cite{CP99}, although there is currently no proof of this property \cite{G03}.
%This delicate behaviour of polygonal billiards with elastic collisions is, however, wiped out as soon as the collisions are non-elastic.
% 
% 

The specific non-elastic rule that we study \cite{AMS09}  shrinks the outgoing angle at a given collision by a uniform factor $0 \le \lambda \le 1$ with respect to the incoming angle at that collision, measured from the normal vector to the boundary. 
 This gives a one-parameter family of maps interpolating between 
 %The parameter $\lambda$ takes values between $0$, where the dynamics reduces 
a one-dimensional map for $\lambda = 0$ and  elastic billiard dynamics for $\lambda = 1$.

% 
% In this case we could naively expect that due to the phase-space contraction, all orbits could end
We might expect that the resulting phase-space contraction would lead only to asymptotic periodic orbits. %  up asymptotic to some periodic orbits.
% It turns out, however, that for triangular billiards nothing could be further from the truth:
In fact, however,   there is also a strong
expansion mechanism present, and the interplay between these  gives rise to hyperbolic (chaotic) strange attractors in phase space.%, with one positive Lyapunov exponent.

This possibility was foreseen in a result of one of the present authors, together with Pujals and Sambarino \cite{MPS10}, that
all pinbilliards on polygonal tables have \emph{dominated splitting}, which is  an extension of hyperbolicity. Pujals and Sambarino \cite{PS09} gave a trichotomy for invariant sets with dominated splitting on a compact manifold: they must be either
finitely many intervals of periodic orbits with bounded period; finitely many simple closed curve where the dynamics is conjugate to an irrational rotation on a circle; or a finite union of hyperbolic sets.

In polygonal tables, the pinball billiard map  is only piecewise continuous,
and hence the phase space is non-compact. For the equilateral triangle, 
we find a non-compact attracting set for any $\lambda$.
%we rigorously show the existence of an attracting hyperbolic set in a certain interval of values of the contraction parameter $\lambda$.
%polygonal pinbilliards two of these phenomena occur: intervals of periodic points (of bounded period) or finitely many hyperbolic invariant sets (some of them can be isolated hyperbolic periodic orbits).
% 
% periodic orbits; conjugate to an irrational rotation on a circle; or a union of hyperbolic sets.  The original motivation of \cite{AMS09} and the present paper
% was to establish which situation pertains in a given case. It turns out, however, that the resulting limit sets are interesting in their own right.
%
Our objective is to study the structure and dynamics of these attractors. The results of the present paper suggest that for any $0 < \lambda < 1$, the attractors exhibit a rich fractal geometry. 
In fact, for $\lambda < \textfrac{1}{3}$ we rigorously prove the existence of a transitive attractor whose structure is roughly a Cantor set of lines, and give a symbolic model that completely describes the dynamics. Moreover,
the attractors for all values of $\lambda$ in this regime are topologically conjugate to one another.
%a strange attractor, i.e.\ an attractor with a non-integer dimension on which the dynamics is chaotic, having a positive Lyapunov exponent (exponentially-fast
 %separation of nearby initial conditions).  
%For any given billiard table, there is thus a  of strange attractors in phase space when $\lambda$ is varied.
%  which is the object of study of the present paper. 

For larger values of $\lambda$, the attractors undergo a sequence of structural changes, giving rise to striking images, shown in fig.~\ref{fig:phenomenology}:
Some regions of phase space become inaccessible due to the nature of the billiard dynamics, which is reflected in the geometry of the attractor. 
A hyperbolic periodic orbit then becomes isolated from the transitive strange attractor, and finally, for $\lambda$ close to $1$, the strange attractor splits into  three distinct transitive components.

\section{Definitions and main results}
% \subsection{Definitions}
\label{sec:defs}

Let $\gamma$ be the boundary of a billiard table, taken to be a closed, connected, convex curve in two dimensions, with length $|\gamma|$, and certain smoothness conditions
%This will be the boundary of a billiard table.
%The required conditions on the smoothness of this curve are standard, as detailed in 
\cite{AMS09, MPS10}.
The boundary of the table is parametrised by arc length $r$.
%  and the billiard dynamics will be given in terms of the outgoing angle $\theta$ at each collision.

\begin{figure}[tp]

\begin{minipage}[c]{0.23\textwidth}
\subfigure[]{
  \includegraphics*[scale=1]{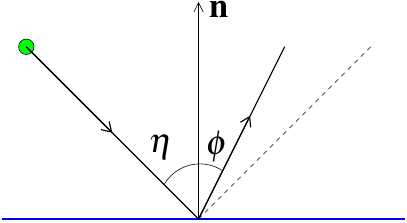}
}
\end{minipage}
\begin{minipage}[c]{0.23\textwidth}
\subfigure[]{
  \includegraphics*[scale=0.38]{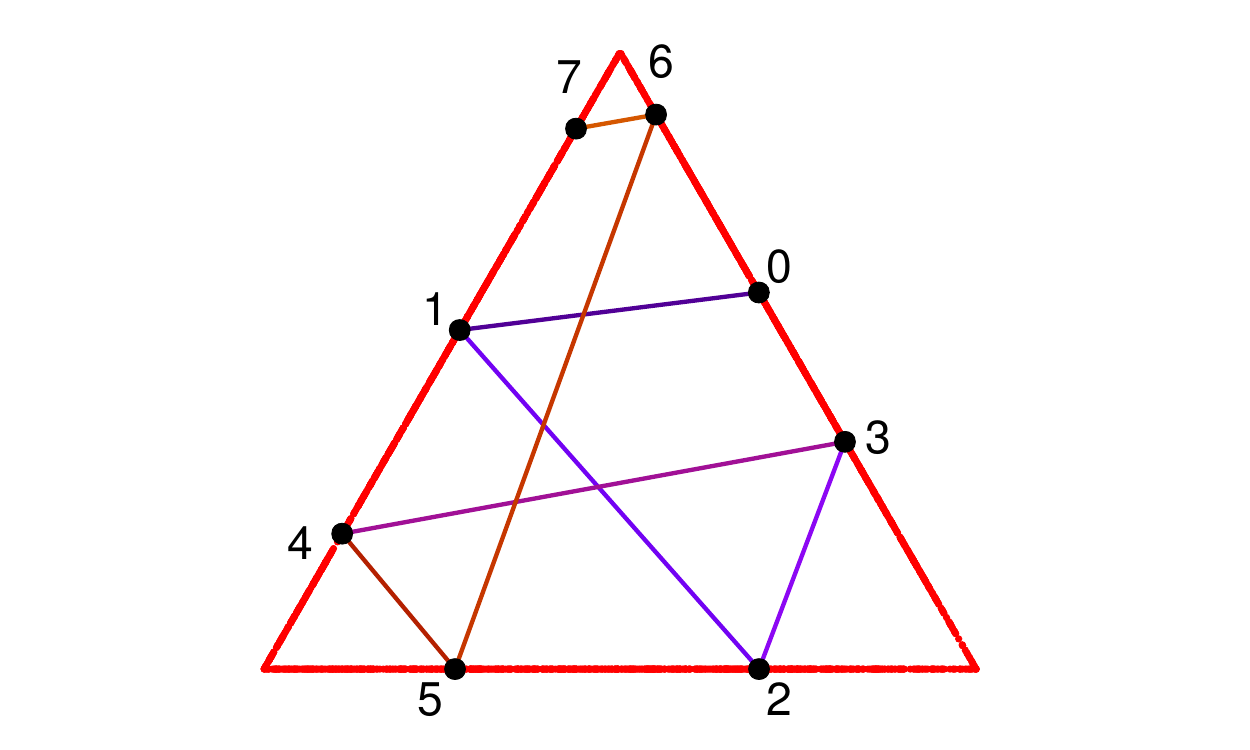}
}
\end{minipage}
\caption{(a) Non-elastic reflection rule;
(b) several bounces in an equilateral triangle for $\lambda = 0.5$, with collisions numbered in order starting from the initial position $0$.
%AGREGAR L y R PARA RIGHT AND LEFT BOUNCE ?
\label{fig:bounce-rule}
}
\end{figure}

We consider non-elastic billiards in which the angle of reflection is modified from an elastic one. General modifications
were previously studied \cite{MPS10}; in  this paper we restrict to the following
one. The
\emph{pinball billiard (pinbilliard) map} $T_\lambda$ with parameter $\lambda$ is given by $T_\lambda(r_0, \theta_0) = (r_1, \theta_1)$,
where $r_1$ is obtained as in the elastic billiard, by moving along the
direction determined by the outgoing angle $\theta_0$, starting at the boundary point $r_0$, until the next boundary collision.  
The new outgoing angle $\theta_1$ (at the next collision) is then given, following a standard sign convention, by $\theta_1 \defeq -\lambda\eta_1$, where $\eta_1$ is
the angle between the incident velocity vector at $q_1$ and the outward normal
$-n(q_1)$; it satisfies
$-\frac{\lambda \pi}{2} \leq \theta_1 \leq \frac{\lambda \pi}{2}$.  This is depicted in figure~\ref{fig:bounce-rule}, together with an example trajectory in an equilateral triangle.

The phase space of the pinbilliard map is thus contained in $\tM \defeq [0, |\gamma|) \times (-\frac{\pi}{2}, \frac{\pi}{2})$.
In fact, there is a finite number of curves where the map is not defined, consisting of points $(s, \theta)$ located at a vertex of the polygon, or whose trajectories hit a vertex in one iteration, so that an outgoing angle cannot be uniquely defined. The singular locus $\DiscontinuityRegion \subseteq \tM$ of $T_{\lambda}$ consists of all preimages of these curves, for which  the dynamics can no longer be continued after some iterate at which the dynamics hits a vertex.
We define the phase space $M \defeq \tM \setminus \DiscontinuityRegion$, the set of points which can be iterated arbitrarily many times. Note that $\DiscontinuityRegion$ has zero Lebesgue measure, since it is the union of countably many regular curves.

The main object of interest is the structure of the maximal forward invariant   $\Gamma_\lambda$ set in $M$, which we refer to as the \emph{attractor} of $T_\lambda$:
\begin{equation}
  \Gamma_\lambda \defeq \bigcap_{n \ge 0} T_\lambda^n(M).
\end{equation}
Note that $\Gamma_\lambda \subset M = \tM \setminus \DiscontinuityRegion$, and  is thus not necessarily compact.
% This is why invariant sets are not compact.???

Except in section \ref{sec:proofMainTheorem}, we shall generally omit reference to the singularity set $\DiscontinuityRegion$, but it is always understood that this set must be excluded when discussing the dynamics.

\subsection{Structural changes of strange attractors}
\begin{figure*}
\subfigure[$\lambda=0.3$]{%
\includegraphics*[scale=.16]{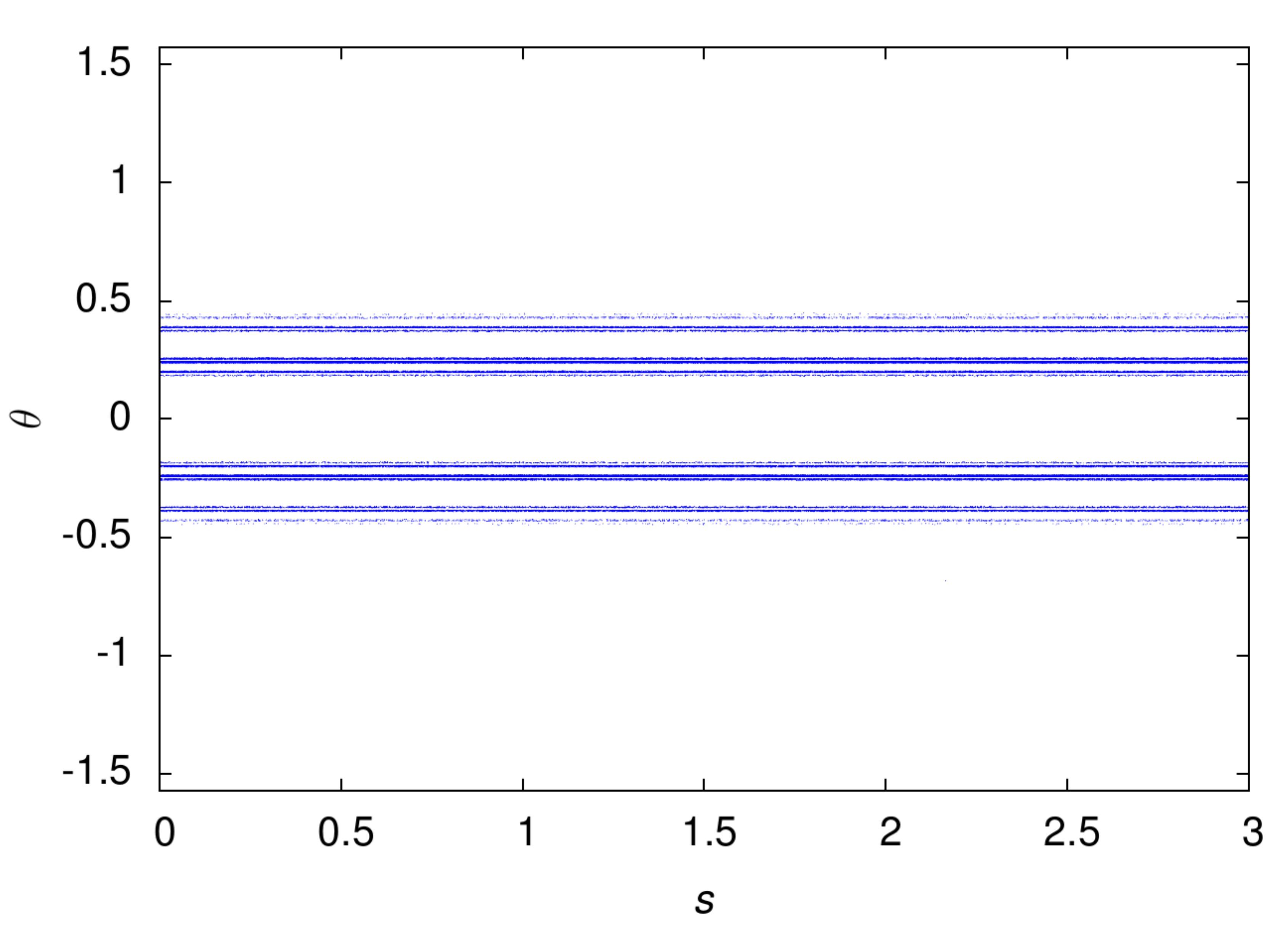} \label{fig:triangle:a}
}
\subfigure[$\lambda=0.5$]{%
\includegraphics*[scale=.16]{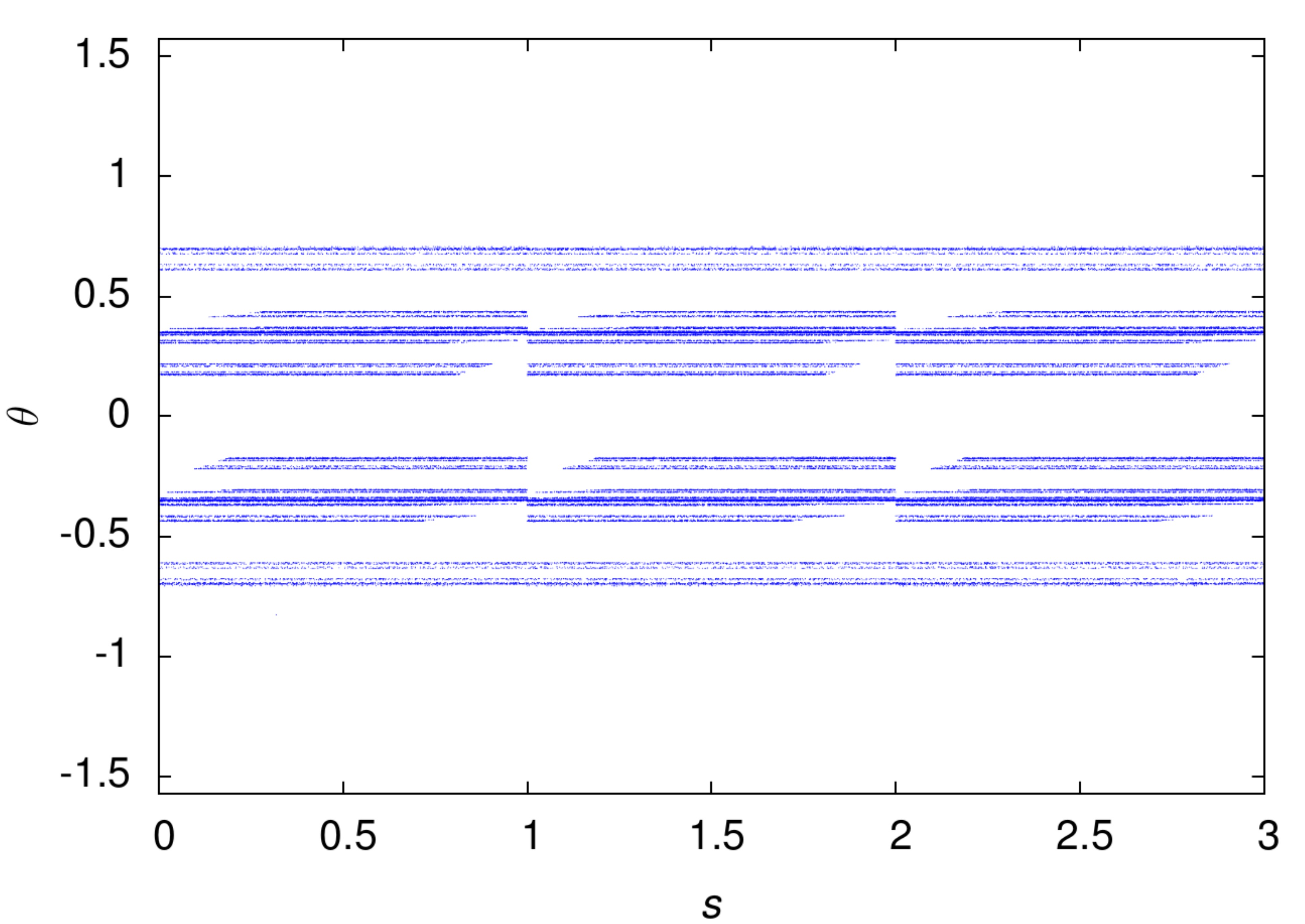}
}
\subfigure[$\lambda=0.68$]{%
\includegraphics*[scale=.16]{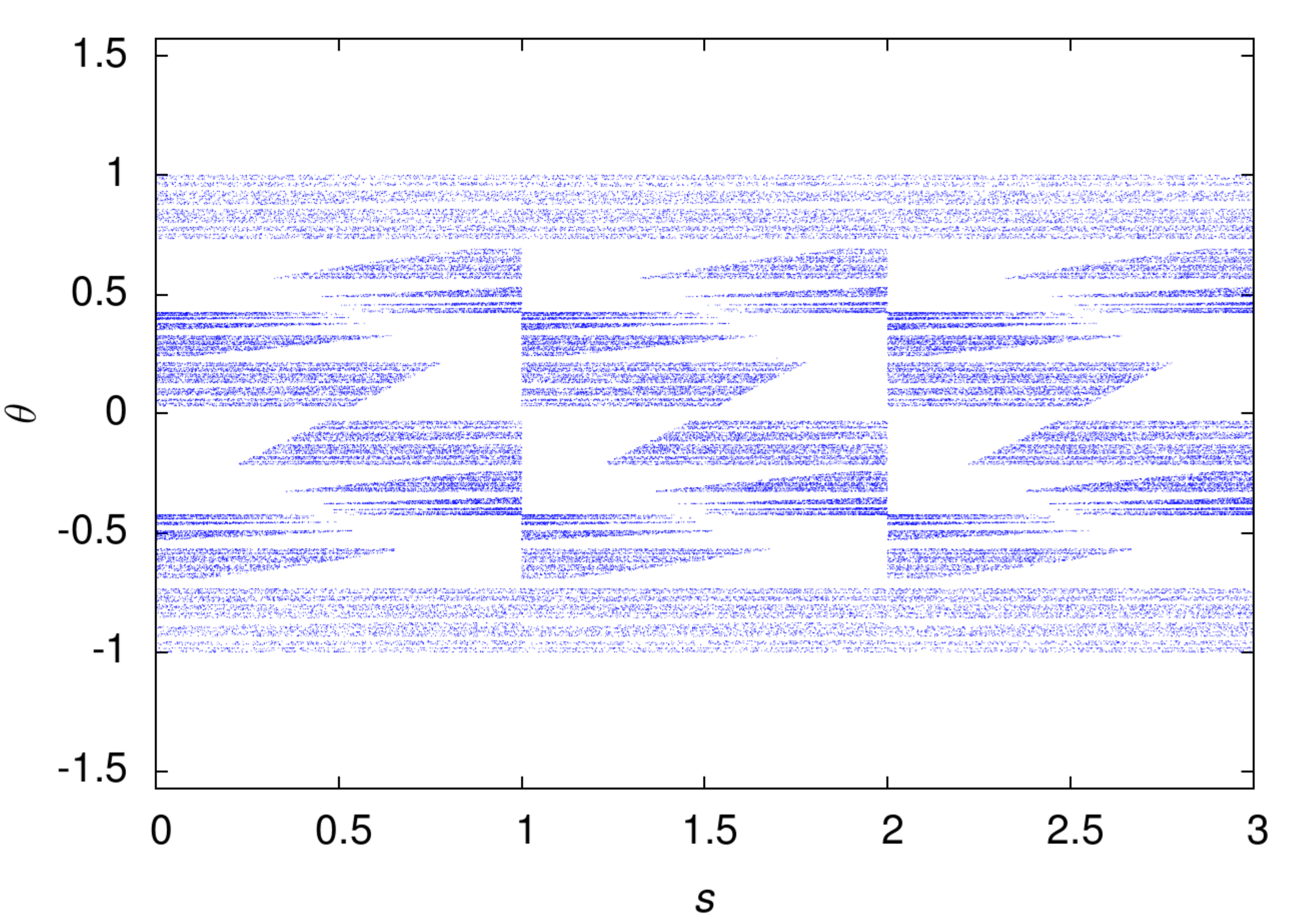}
}

\subfigure[$\lambda=0.8$]{%
\includegraphics*[scale=.16]{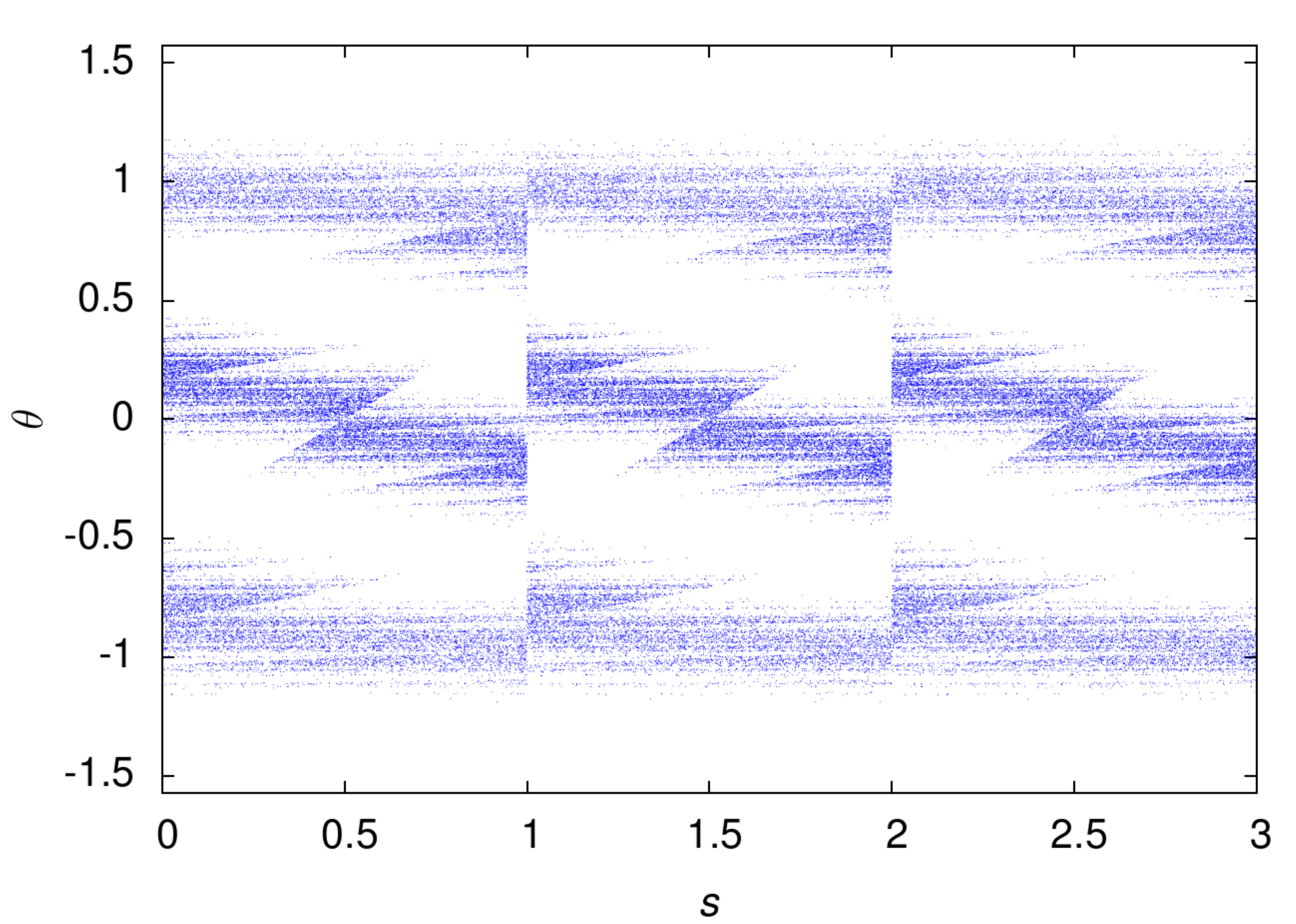}
}
% \qquad
\subfigure[$\lambda=0.95$]{%
\includegraphics*[scale=.16]{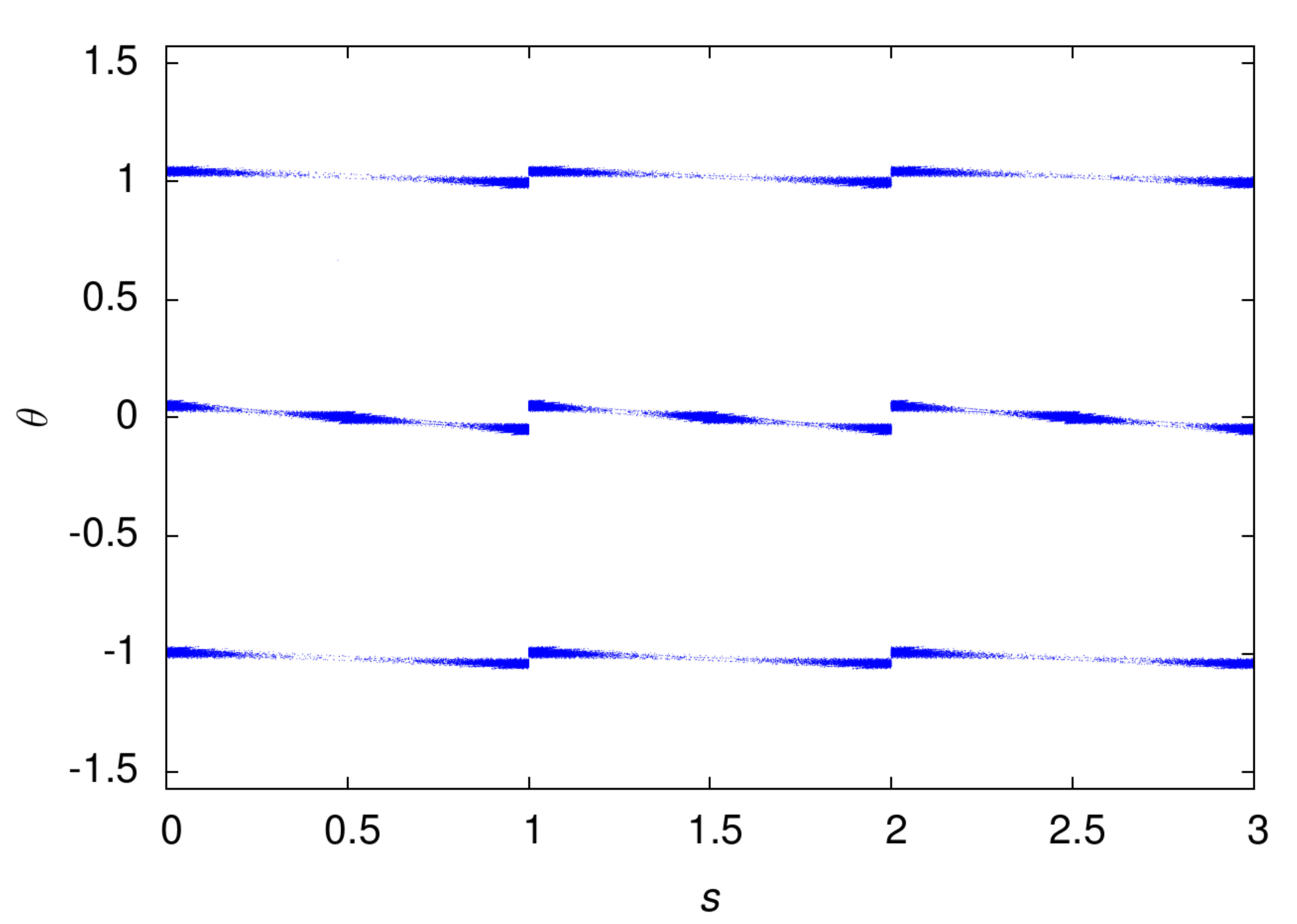}
}
%\quad
\subfigure[$\lambda=0.98$]{%
\includegraphics*[scale=.15]{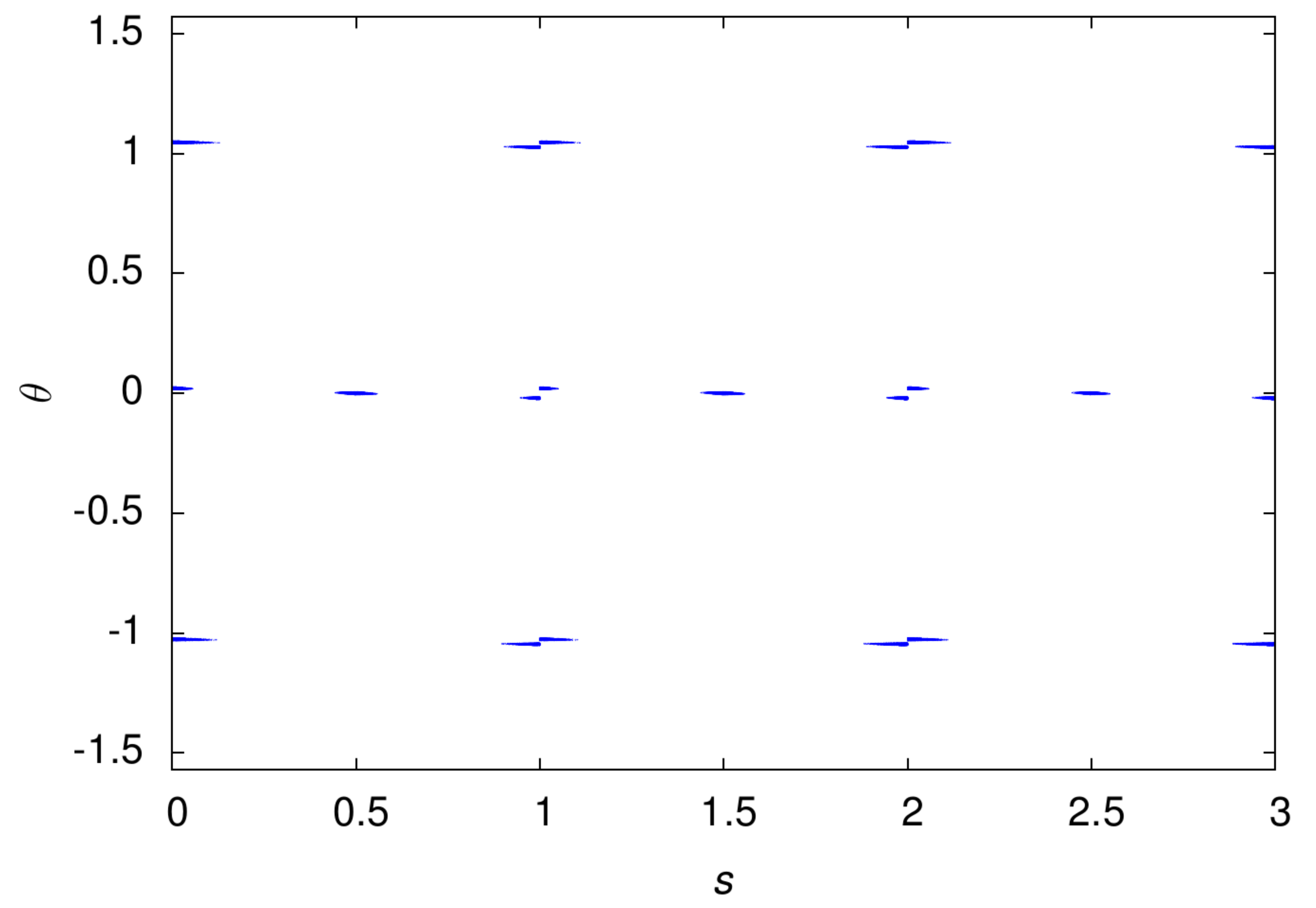}
}

\caption{Phase space of pinbilliard dynamics in the equilateral triangle, depicting the main structural
changes observed in the strange attractor as the contraction parameter $\lambda$ is varied.
Each figure except (f) shows a single orbit of length $10^5$ starting from a random initial condition, after a transient of
 $10^{4}$ iterates were discarded to reach the asymptotic attractor. In (f), three such orbits are shown. (Enhanced online.)
}
\label{fig:phenomenology}
\end{figure*}

Figure~\ref{fig:phenomenology} shows numerically-obtained attractors in phase spacefor pinbilliard dynamics in an equilateral triangle of side length $1$, for several values of $\lambda$. We see that trajectories from almost all initial conditions are observed to accumulate on fractal strange attractors. The main goal of this paper is the explanation of this phenomenology.
The main structural changes observed are the following.
% 
% depicting the following main structural changes.
% For pinbilliard dynamics in the equilateral triangle, trajectories from almost all initial conditions are observed to accumulate on fractal strange attractors.
% % as $\lambda$ varies between $0$ and $1$.
% In this section, we exhibit the observed structural changes of the attractors found for pinbilliard dynamics in the equilateral triangle
% as $\lambda$ varies between $0$ and $1$.

% DRAW ON 3-periodic orbit

% Figure~\ref{fig:phenomenology} shows numerically-obtained phase portraits for pinbilliard dynamics in the equilateral triangle for several values of $\lambda$, depicting the following main structural changes.

For sufficiently small $\lambda$, in fact $\lambda < \frac{1}{3}$, the attractor is a Cantor set of lines, which ``thickens'' as $\lambda$ increases.
For larger values of $\lambda$,  triangular-shaped gaps begin to appear in the attractor (b), which become larger as $\lambda$ increases.
Afterwards, around $\lambda \simeq 0.68$, the complete horizontal lines visible in the central bands of the attractor break up (c), which we refer to as band splitting. 
Soon thereafter, round $\lambda \simeq 0.7$, there is band merging to form just three bands (d). These bands shrink as $\lambda$ increases, concentrating around three angle values (e). Finally, for $\lambda \simeq 0.98$ and above, the attractor breaks into three distinct transitive regions (f). The enhanced online multimedia version of this figure shows the evolution of the attractor for all values of $\lambda$.
These phenomena will be analyzed in turn in subsequent sections.
%, separately, from section \ref{sec:equiTriangle} to section \ref{sec:nontransitive}.
% and in section \ref{sec:proofMainTheorem} a proof of Theorem \ref{MainTheorem} will be given.

% \section{Phenomenology of strange attractor evolution}

% Numerical observations -- see figure~\ref{fig:phenomenology} -- exhibit a rich structure for the attractors in this case, with distinct types of behaviour in different subintervals of $\lambda$.
\subsection{Main theorem for $\lambda < \textfrac{1}{3}$}
For $\lambda < \textfrac{1}{3}$, we are able to characterise the dynamics completely via a symbolic model:

\

% Based on these observations, we find the following main results.
% There turn out to be several different types of behaviour in different sub-intervals of $\lambda$. 
% 
% 
% 
% For small enough $\lambda$, in fact $\lambda \in (0, \frac{1}{3})$, the
% attractor is shown by the following theorem to be roughly a Cantor set of lines:
% , and the dynamics is a non-linear
% version of a ``skinny baker map'' REF. %-- see Theorem 1 in Section~???
% In fact, the set of possible angles in the attractor is related to the structure
%  of random geometrical series, as given by the following % (infinite Bernoulli convolutions):

% for a given $0< \lambda <1$, consider the set $C(\lambda) \defeq \{ \sum_{n\geq 1} \pm \lambda^n\} \subseteq \RR$.
% For $0<\lambda < \frac{1}{2}$, this is a Cantor set of zero Lebesgue measure in $\RR$.
% 
% %In a certain region of the parameter $\lambda$, w
% We have the following strong
% 
% 

\begin{theorem}[Main Theorem] \label{MainTheorem}
Let $0< \lambda < \textfrac{1}{3}$. Then the attractor $\limitset$ of the pinball billiard map  in the equilateral triangle is a transitive invariant set and is homeomorphic to $[0,3] \times C \setminus \DiscontinuityRegion$, where $C$ is a Cantor set in $\RR$.
%Moreover, it is hyperbolic, and the horizontal lines 
Moreover, there is an invariant ergodic measure on $\limitset$, which is absolutely continuous with respect to the Lebesgue measure on horizontal lines.
\end{theorem}

Section~\ref{sec:proofMainTheorem} of this paper is devoted to a rigorous proof of the Main Theorem.

%\noindent \emph{Remark:} 
%For $\lambda \in [0,\textfrac{1}{3}]$, the set $\limitset$ is  a hyperbolic invariant set of %saddle type. 

%Recall that a one-parameter family of dynamical systems is structurally stable in an %interval of the parameter if they are mutually topologically conjugate.

Note that since all of these maps, restricted to the attractor, are conjugate to the same model, they are also conjugate to one another, i.e.\ the one-parameter family of pinbilliards is structurally stable for $\lambda<\frac{1}{3}$.

\section{Equilateral triangle} \label{sec:equiTriangle}

%In this section we begin to analyse the details of pinbilliard dynamics in the equilateral triangle.
%

%In this paper we study pinbilliard dynamics in an equilateral triangle.
For the pinbilliard map that we will study, we consider an equilateral triangle of
side length $1$. %, so that the arc length $s$ lies in the interval $[0,3]$.
The possible phase space coordinates are then $\tM = (0,3) \times (-\textfrac{\pi}{2}, \textfrac{\pi}{2})$, which we split into $\cup _{i=0}^2 M_i$, where $M_i \defeq (i, i+1) \times (-\textfrac{\pi}{2}, \textfrac{\pi}{2})$ is the part of phase space corresponding to the side $i$, 
with arc length in $(i, i+1)$.  (Note that the dynamics is not defined for phase space points whose arc-length parameter sits on a vertex, so these points are automatically excluded from the phase space.)
%The pinbilliard map $T_{\lambda}:M\to M$.
%The limit set of $T_{\lambda}$ will be denoted by $L_{\lambda}$.
When $i$ denotes a side of the triangle, it must be read modulo $3$, i.e.\ $i\in \ZZ_{3}$. %; e.g.\ $i= -1 \equiv 2 \mod 3$ refers to the third side, 
%starting from $0$. 
%In the case of the equilateral triangle, the phase space must inherit its rotational symmetry.
In these coordinates, we will refer to vertical and horizontal 
intervals as being lines in phase space with constant first or second coordinate, respectively.

\subsection{Discontinuity curves} \label{subsec:discontinuity}

\begin{figure}[tp] %  figure placement: here, top, bottom, or page
   \centering
   \includegraphics[scale=0.4]{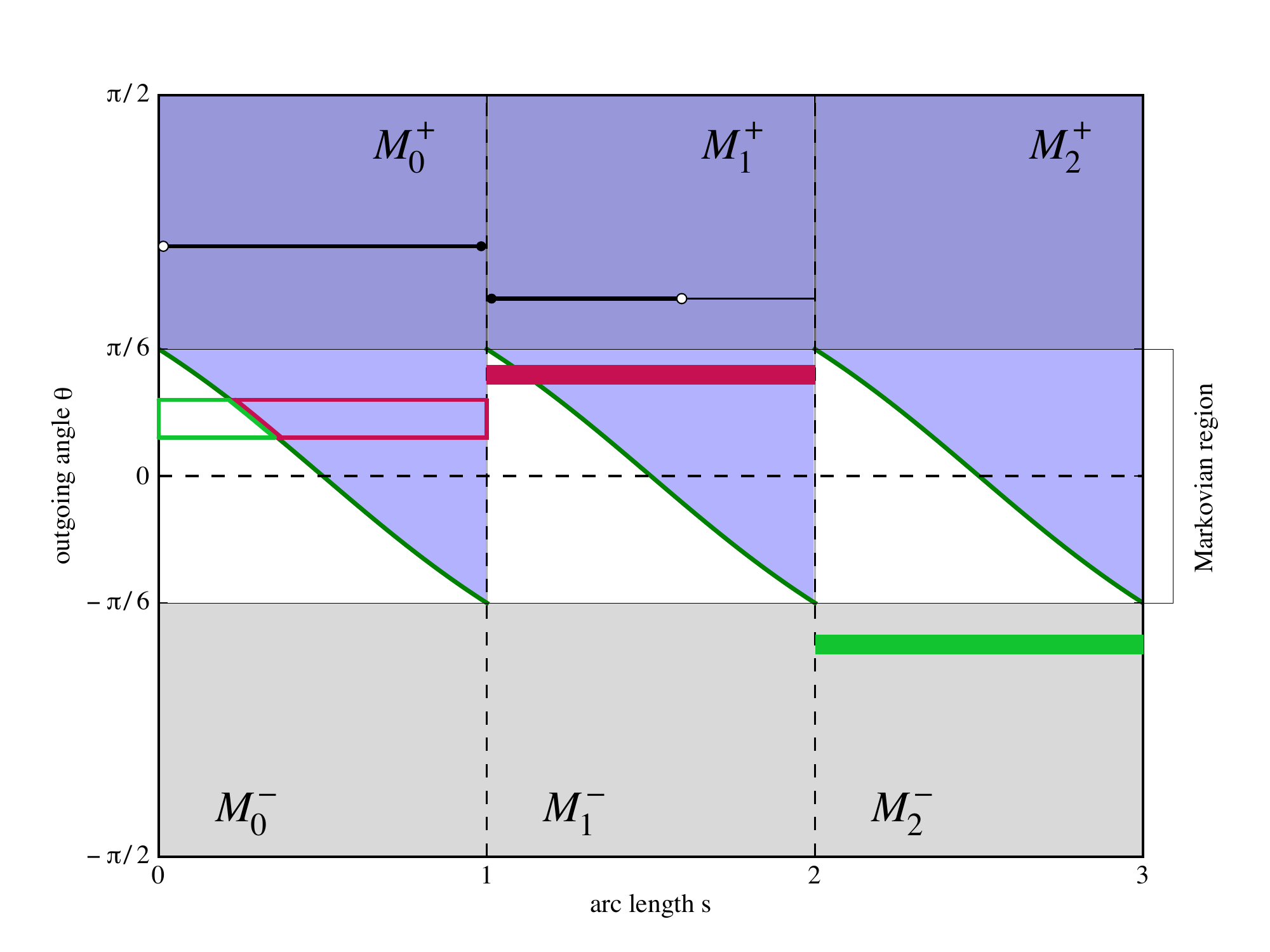}
   \caption{%
	Schematic picture of the phase space of the equilateral triangle, showing the discontinuity lines (solid curved lines) and the resulting domains of continuity $M^{\pm}_i$ (shaded and empty).
	The Markovian region is shown, in addition to a rectangle (green and red open boxes) together with its images (green and red solid rectangles). The horizontal line shows the non-Markovian dynamics, with the open and closed circle indicating the orientation-reversing feature.
The dashed vertical lines indicate the positions of the vertices. 
   }
   \label{fig:dominios}
\end{figure}

Key objects in billiards with piecewise-smooth boundary are the \emph{discontinuity curves}, %referred to above in Section~???, 
which consist of those phase space points that are mapped to a vertex.
For the equilateral triangle, the discontinuity curve consists of points $(s, \delta(s)) \in M$ with 
\begin{equation} \label{eq:curvasdiscontinuidad}
\delta(s) \defeq \arctan \left(  \frac{1+2(\lfloor s \rfloor - s)}{\sqrt{3}}   \right), %\quad i= 0, 1, 2
\end{equation}
where $\lfloor s \rfloor$ denotes the integer part of $s \ge 0$. 
The function $\delta$ is 
independent of $\lambda$, and has image
$(-\frac{\pi}{6}, \frac{\pi}{6})$. 
As we shall see, this curve splits each $M_i$ into two parts
$M^+_{i}$ and  $M^-_{i}$, which are the domains of continuity of $T_\lambda$; see figure~\ref{fig:dominios}.

Since the geometry of the equilateral triangle is sufficiently simple, we can compute explicitly the angular coordinate of $T_\lambda$ in terms of two affine maps:
%
%The importance of this curve stems from the fact that it allows us to compute explicitly the action of the pinbilliard map on the  the angular (second) coordinate of the pinbilliard map then takes the following explicit form:
\begin{equation}T_{\lambda}(s, \theta) =
\left \{
\begin{array}{l}
( \,\cdot \, , \varphi_{+1}(\theta))  \qquad \text{ if } \quad \theta > \delta(s),  \\ 
( \,\cdot \, , \varphi_{-1}(\theta))  \qquad \text { if } \quad \theta < \delta(s),
\end{array}\right.
\label{eq:T-explicit}
\end{equation}
where
the functions $\varphi_{\pm 1} : \RR \to \RR$ are defined by % for $z \in \{+1, -1\}$ by
\begin{eqnarray}
\varphi_{\pm 1}(\theta) 
\defeq \lambda \left( \pm \textstyle \frac{\pi}{3} - \theta \right).	\label{phi+-}
\end{eqnarray}
%
%as given by the pinball law of collision.  
To see this, note that from a starting point $(s, \theta) \in M$, 
%To see this, note that from a given starting position $s$,
 on the side $i  \defeq \lfloor s \rfloor$ of the triangle, there are only two possible types of bounce, depending on $\theta$: a \emph{right bounce}, on side $i+1$, or a \emph{left bounce} on side $i-1$; these possibilities are separated by the point $(s, \delta(s))$  that hits the vertex.  A straightforward computation verifies that the angular coordinate after the next collision is given by $\phi_z(\theta)$, where the sign $z \in \{+1, -1\}$  denotes a right or left bounce, according to the sign of the angle $\pm \frac{\pi}{3}$ between the sides  before and after the bounce.

% 
% For any angle $\theta$ in the $(-\frac{\pi}{6}, \frac{\pi}{6})$, both left and right bounces can thus be obtained starting from the initial condition $(s, \theta)$ for distinct intervals of points $s$, and 
%  the corresponding angular coordinate of $T_{\lambda}(s, \theta)$ is $\varphi_{\pm 1}(\theta)$. 

The form of $T_\lambda$ in \eqref{eq:T-explicit} implies that locally, horizontal intervals  in phase space  are mapped onto other horizontal intervals with in general different angles, given by projections along a fixed direction 
from one side to another.

\subsection{Derivative}

To calculate expansion and contraction properties of the pin billiard map $T_{\lambda}$, it is necessary to calculate its derivative at a point $(s_0, \theta_0)$. 
The expression for this obtained for a general geometry\cite{MPS10} reduces for polygonal billiard tables, which have only flat boundary components,  to
\begin{equation} \label{mpders} 
D_{(s_0, \theta_0) }T_\lambda = - \left(
\begin{array} {ccc} 
A & B
\\ 0   & \lambda
\end{array} \right), \quad 
 \text{ with} \quad  A = \frac{\cos \theta_0}{\cos\eta_1},  \quad B = \frac{t_0}{\cos \eta_1};
 \end{equation}
$t_{0}$ denotes the distance between bounces in the billiard table and $\eta_1$ the incoming angle at the next collision, as in sec.~\ref{sec:defs}.

The upper-triangular structure of the derivative of $T_\lambda$ implies that 
 horizontal tangent vectors of the form $v = (a, 0)$ are mapped back into horizontal tangent vectors with different base point.

\subsection{Periodic orbits}
\label{sec:periodic-orbits}

 %This allows us to conclude hyperbolic properties for periodic orbits.
Billiard maps do not have fixed points, and moreover in polygonal billiards neither are there period-$2$ orbits, unless two sides are parallel, in which case there is in fact a whole interval of such orbits. For pinbilliards with our rule, this still holds; thus, any periodic orbit in the equilateral triangle must have period at least $3$. 

\

%
%$ M \subset R^2$, as a differentiable manifold, the tangent bundle is $TM = M \times R^2$, the full derivative map $DT: TM \to TM$ is given by $(T_{\lambda}(s,\theta), D_{(s,\theta)}T(v))$. 
%
%\begin{lemma}
%The linear bundle given by horizontal subspaces on each point of $M$, is invariant under $DT$.
%%The derivative leaves invariant the horizontal sub-bundle of the tangent bundle of M. 
%\end{lemma}

%%The image of the horizontal tangent vector $(1,0)$ by $D_{(s_0, \theta_0)}T_\lambda ^n$ is then another horizontal vector, with first coordinate equal to

 \begin{lemma}  
 Fix $0 < \lambda < 1$. Then
any periodic orbit of $T_\lambda$ is hyperbolic of saddle type.
%has a horizontal expanding direction.
%Furthermore, it is hyperbolic of saddle type if $\lambda < 0.65$. 
%  
% In polygonal pinbilliards, the lines in phase space of constant angles corresponding to
% $n$-periodic points with $n>2$ are invariant and uniformly expanded by $DT^n_\lambda$. 
\end{lemma}

\proof
Consider a periodic point $(s_0, \theta_0)$ with period $n \ge 3$ and orbit
$(s_j, \theta_j) \defeq T_\lambda^j(s_0, \theta_0)$ for $j=1, \ldots, n$.
The periodicity implies that $\theta_n = \theta_0$.

The  derivative $D_{(s_0, \theta_0)} T_\lambda^n$ is an upper-triangular matrix 
with eigenvector $(1,0)$. The absolute value of the corresponding eigenvalue is
%whose action on the horizontal tangent vector $(1,0)$ by multiplying it by
\begin{equation} \label{der-uns}
\mu_n \defeq \frac{\cos \theta_0}{\cos \eta_n} \prod _{i=1}^{n-1}
\frac{\cos \theta_i}{\cos \eta_i} =
\frac{\cos \theta_0}{\cos \eta_n} \prod _{i=1}^{n-1}
\frac{\cos \lambda \eta_i}{\cos \eta_i}, \end{equation}
where $\eta_i$ is the incoming angle at the $i$th collision.
The left-hand side follows from the chain rule and eq. \eqref{mpders} after a rearrangement; the equality follows from the pinbilliard collision rule $\theta_i = -\lambda \eta_i$, and using the fact that $\cos$ is an even function.
Finally, we have $\cos \theta_{0} / \cos{\eta_{n}} = \cos (\lambda \eta_{n}) / \cos(\eta_{n})$, so that $\mu_{n} = \prod_{i=1}^{n} [\cos(\lambda \eta_{i}) / \cos \eta_{i}] > 1$.
Thus the line spanned by $(1,0)$ is the unstable eigenspace of the periodic orbit.

%
%Considering horizontal lines in phase space with constant angle, 
%we obtain the following
% \begin{lemma}  In polygonal pinbilliards, the lines in phase space of constant angles corresponding to
% $n$-periodic points with $n>2$ are invariant and uniformly expanded by $DT^n_\lambda$. 
%\end{lemma}
%

Since the product of matrices of this form along the periodic orbit is also upper-triangular,
its determinant is equal to the product of its diagonal entries and is also equal to the product of its eigenvalues.  The other eigenvalue is thus $\lambda^{n} < 1$. Hence, there is one expanding and one contracting direction, and thus the periodic orbit is hyperbolic of saddle type.
%
%Since $\theta_0 = \theta_n = \lambda \eta_n$, each term in \eqref{der-uns} is greater than $1$, since $\cos$ is a decreasing function and $|\lambda \eta_i| < |\eta_i|$ when $|\eta_i| < \frac{\pi}{2}$.  
%
%To prove that such a periodic orbit is hyperbolic of saddle type, we must show that it has another eigenvalue that is less than $1$ in absolute value; a sufficient condition for this is that  the determinant has absolute value less than $1$.
%We have $|\det D T_\lambda^n| = \lambda^n \mu_n$. If $\lambda \cos (\lambda \eta)  / \cos \eta < 1$ for all $\eta \in [-\lambda \frac{\pi}{2}, \lambda \frac{\pi}{2}]$, then each term in this product is less than $1$, and hence the determinant also is. By direct inspection of this function, we find that this is satisfied for any $\eta$ only under the condition $\lambda < 0.65$.  
\qed

Note that this result holds for any periodic orbit of any polygonal pinball billiard,
except the period-$2$ orbits mentioned above.
%For larger values of $\lambda$, detailed knowledge of the complete orbit must be taken into account to guarantee hyperbolicity.

%
%\proof The invariance follows from the fact that lines of constant angle map under
%$T_\lambda$ to lines of constant angles. For periodic points, 
% \eqref{der-uns} becomes $ \prod _{i=1}^{n}
%\frac{\cos \lambda \eta_i}{\cos \eta_i} > a^n,$ for some $a > 1;$ this implies
%the uniform expansion.
%\qed
%\begin{remark}
%This last lemma holds for any polygonal billiard table.
%\end{remark}

\subsection{Distinguished 3-periodic orbits}\label{3per}
For \emph{elastic} billiard dynamics in an equilateral triangle, there are two distinguished period-$3$ orbits,
formed by the median triangle and an orientation. It turns out that we may continue these period-$3$ orbits for pinbilliard dynamics for all  $\lambda < 1$. They maintain their rotational symmetry and we have the following lemma:
%he continuation of this period-$3$ orbit is also present: 
%$\lambda < 1$ always exists  pinbilliards for all
%$\lambda \in [0,1]$:
\begin{lemma}
There are two hyperbolic period-$3$ orbits of saddle type for any $\lambda \in (0,1)$. Their local unstable manifolds are horizontal lines.
% contained in horizontal lines.
\end{lemma}

\proof 

The rotational symmetry of the distinguished period-$3$ orbit suggests that the outgoing angle must be a fixed point $\theta_*$ of the angular dynamics, so that  
$\varphi_{+ 1}( \theta_{*}) = \theta_{*}$.
It is a straightforward computation that the unique solution is $\theta_* \defeq \theta_*(\lambda) \defeq \frac{\pi}{3} \frac{\lambda}{\lambda +1}= \frac{\pi}{3}(\lambda - \lambda^2 + \lambda^3 - \cdots)$. 
In fact,  $\varphi_{-1}(- \theta_{*}) = -\theta_{*}$ is the
counterpart with reversed orientation.
%Note that $0 < \theta_* < \frac{\pi}{6}$ for any $0 <\lambda < 1$.
%This shows that there is an angle which is invariant under the dynamics.
%
%For there to be a period-$3$ point with arc-length coordinate $s_3(\lambda)$ and angle $\theta_*$, we must 

A simple trigonometric calculation using the rotational symmetry of the period-$3$ orbit shows that the point $(s_3(\lambda), \theta_*(\lambda))$ with
\begin{equation}
%  s_3(\lambda) \defeq \left(\frac{2}{\sqrt{3}\tan(\theta_*) +3 },\theta_* \right) \in M^+_{0}.
 s_3(\lambda) \defeq \frac{2}{\sqrt{3}\tan(\theta_*) +3 },
\end{equation}
is periodic of period $3$ for $T_{\lambda}$.

%Since $|\theta^*| < \frac{\pi}{6} < 0.65$  for all $\lambda \in (0,1)$, 
The previous Lemma then shows that this period-$3$ orbit is always hyperbolic.

To see that the local unstable manifold of the period-$3$ point is a horizontal interval, let $J \subset [0,3] \times \{\theta_*\} $ be a small horizontal interval that contains the period-$3$ point. 
Then $T_\lambda^3 (J) \supset J$, so $J$ is a local unstable manifold of the period-$3$ point.
\qed

%Note that $\pm \theta_{*}$ is an attracting fixed point of $\varphi_{\pm 1}$, since
% $|\varphi_{\pm 1}'(\pm \theta_*)| = \lambda < 1$.

\section{Markovian region and complete symbolic dynamics}
\label{sec:markovian-heuristic}
%We shall exploit the following geometrical feature of polygonal billiards:
%The image of a horizontal interval, i.e.\ one with constant outgoing angle, 

\begin{figure*}[tp] %  figure placement: here, top, bottom, or page
\subfigure[Markovian]{%
\includegraphics*[scale=0.08]{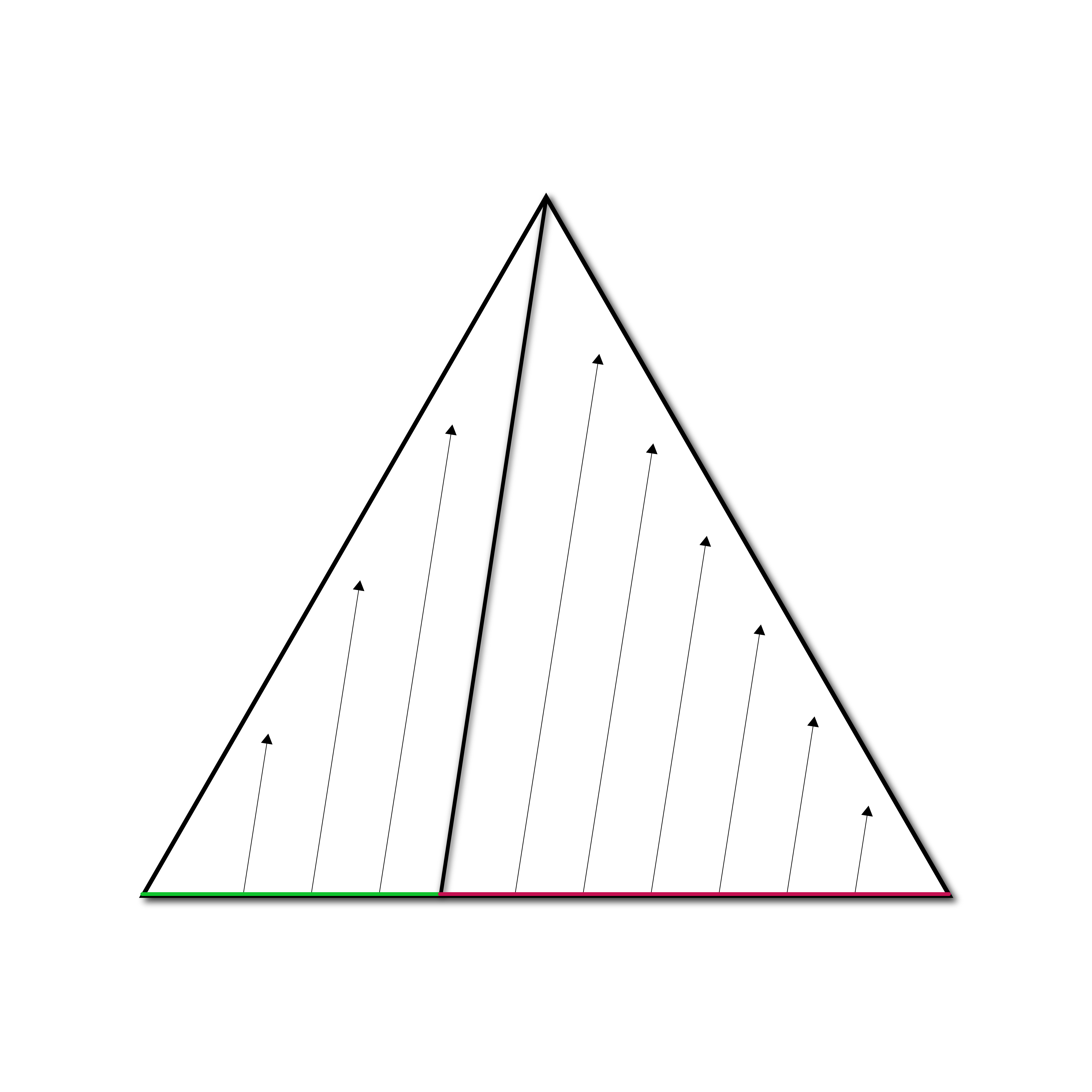} 
}
\quad
\subfigure[Non-Markovian]{%
\includegraphics*[scale=0.08]{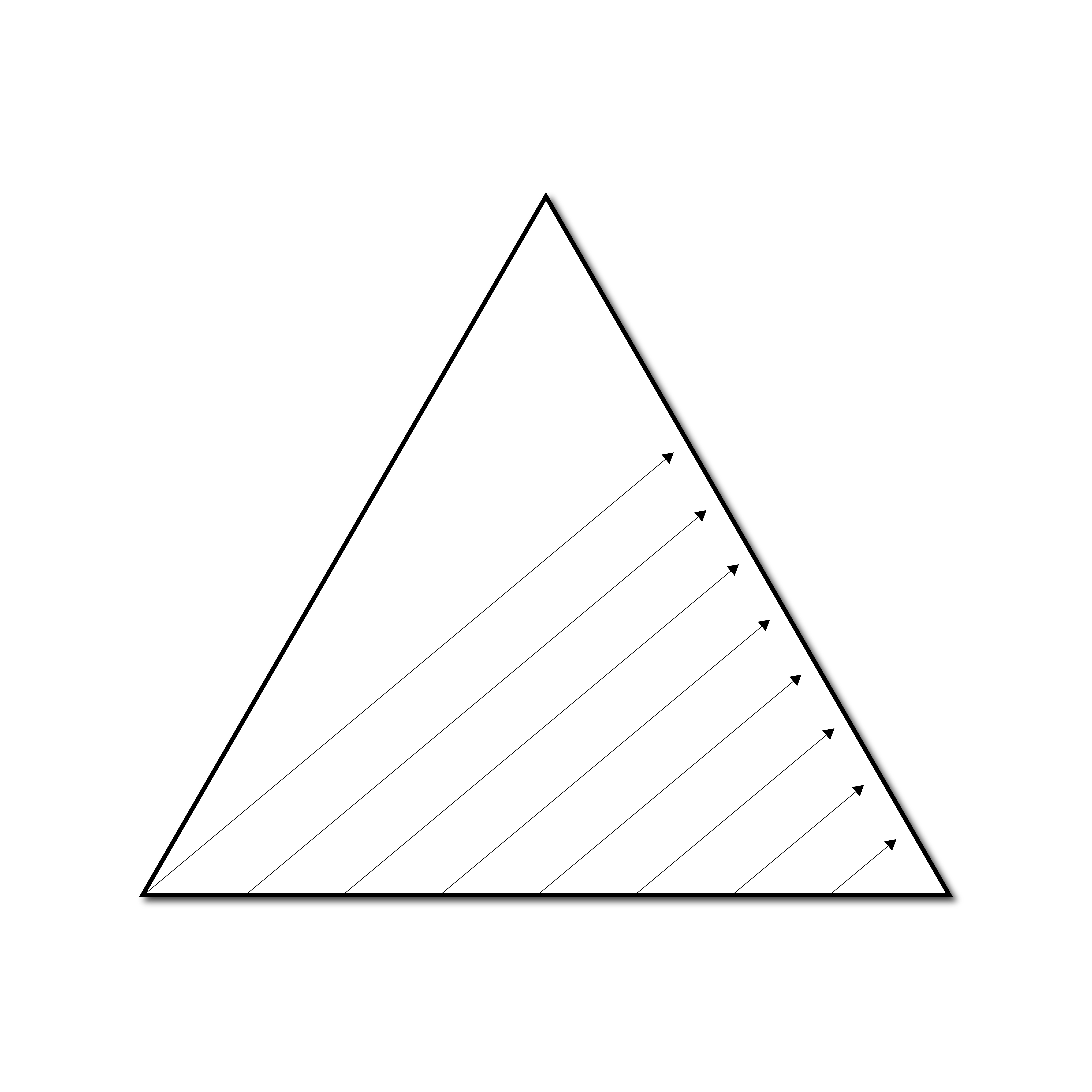}
}
\caption{%
Images of horizontal intervals: (a) Markovian case, in which the other two complete sides are covered and there is linear expansion; (b) non-Markovian case, in which only one partial side is covered, and there is linear contraction. \label{fig:markschematic}
   }
   
\end{figure*}

We shall use the following feature of billiard dynamics in an equilateral triangle:
Let $J_\theta$ be a horizontal interval which covers one side of the triangle:
% to obtain a complete description of the dynamics in the equilateral triangle for sufficiently small 
%$\lambda$ using a symbolic model, 
%Consider any horizontal interval in phase space; that is an interval contained in the boundary of the table and with a fixed angle $\theta \in (-\textfrac{\pi}{2}, \textfrac{\pi}{2})$.  The first iteration of such an interval is given by a finite union of horizontal intervals, which are the images of affine projections between corresponding sides of the polygon, along a fixed direction given by the angular coordinate $\theta$. The angular coordinate of the image (outgoing angle) changes, but depends only on which side of the table the orbit bounces off.
% 
%For the equilateral triangle, this can be seen from eq.~\eqref{themap} ???, since the second coordinate of the pinbilliard map is given by two affine functions which depend only on the angle, and the one which is used depends on the position of the point with respect to the discontinuity curves; see eq.~\eqref{eq:curvasdiscontinuidad}. 
 %
%In general, two kinds of situations happen for horizontal intervals of the form 
$$J_\theta := (i, i+1)\times \{\theta\} \subset M_{i},$$ with $i = 0$, $1$ or $2$. 
%depending on $|\theta|$, but 
There are two possibilities for the structure of the image $T_\lambda(J)$, depending on the value of $\theta$:
either it is a single horizontal interval or a union of two horizontal intervals, depending on whether the projection in the direction $\theta$ of the side $i$ of the triangle does or does not hit a vertex; see fig.~\ref{fig:markschematic}.

When $|\theta| > \textfrac{\pi}{6}$, the image $T_{\lambda}(J)$ is a horizontal interval contained in $M_{j}$, where $j$ corresponds to the previous or next side in the table, depending on the sign of $\theta$; and $T_{\lambda}(J)$ is linearly contracted by a factor that  depends only on $\theta$ (and not on $\lambda$), and which tends to $1$ as $|\theta| \to \textfrac{\pi}{6}$. % in fact, for this limiting case $T_{\lambda}|_{J}:J\to T_{\lambda}(J)$ is an isometry. 

On the other hand, when $|\theta| < \textfrac{\pi}{6}$ the situation is very different.
The interval $J$  is now split by the discontinuity curve $\delta$, and hence its image consists of two horizontal intervals, each of which completely covers its corresponding side $j \pm 1$. Moreover, each is linearly expanded and maps across the whole horizontal extent of one of the other two sides. %Each is mapped onto one of the two other sides and their images are uniformly expanded. 
For this reason we call  the subset $[0,3] \times (-\textfrac{\pi}{6}, \textfrac{\pi}{6}) \subset M$ the ``Markovian region". 
Note that these geometrical properties in the Markovian region in no way depend on $\lambda$.

The key observation is now that when $\lambda < \textfrac{1}{3}$, the image  $T_{\lambda}(M)$ of the whole phase space is in fact contained in the Markovian region, and hence, roughly speaking, we can obtain a good partition in terms of future itineraries -- i.e., whether the bounce is a left or right bounce at each iteration.
%repeating recursively this previous argument. 
This partition is described precisely in section~\ref{markovianReg}.

\subsection{Angular dynamics}

The combinatorics of the horizontal intervals described in the introduction of this section depend on the dynamics of the angular coordinate, which we proceed to study in this section.

The form of the maps $\phi_{\pm}$ in \eqref{phi+-} leads us to consider the following geometric series with random signs. Let  $0<\lambda <1$ and consider the following subset of real numbers:
\begin{equation} \label{ElCantor}
\elcantor \defeq \left\{ \frac{\pi}{3} \displaystyle \sum_{n=1}^\infty z_n \lambda^n :  z_n \in \{1,-1\} \, \forall n \geq 1 \right\}.
\end{equation}
%ROBERTO TIENE DUDAS SOBRE LA SIGUIENTE FRASE 
We show in section \ref{sec:proofMainTheorem} that $C(\lambda)$ is a Cantor set for $\lambda < \frac{1}{2}$.

Suppose that we are given an angle $\theta \in \elcantor$, that is $\theta = \frac{\pi}{3} \sum_{n=1}^\infty z_n \lambda^n $  for some choice of signs $z_{n} \in \{+1,-1\}$. Then $\phi_{\pm}(\theta) \in \elcantor$. In fact, 
\begin{equation}\label{cantorInvariante}
\angmap_{\pm1}( \ang)  =  \frac{\pi}{3} \left( \pm \lambda-    \sum_{n=2}^\infty z_{n-1} \lambda^{n} \right). \end{equation}

Thus the images $\angmap_{\pm1}(\elcantor) \subset C(\lambda)$  for both maps $\angmap_{\pm1}$, and hence
$$T_{\lambda}([0,3]\times \elcantor) \subset [0,3]\times \elcantor,$$
 (except for  points in $\DiscontinuityRegion$, where $T_\lambda$ is not defined). Thus a good candidate for an invariant set is that it be contained in the set $[0,3] \times C(\lambda)$, which for $\lambda < \frac{1}{3}$ is a Cantor set of horizontal lines, and this corresponds with the pictures of the computed attractor shown in fig. \ref{fig:phenomenology}(a). In section \ref{sec:proofMainTheorem} we rigorously prove that the attractor is indeed precisely $[0,3] \times \elcantor$ without the singularity set.
 
%We thus see that the precise statement of the main result for
% $\lambda <\textfrac{1}{3}$, that of Theorem~\ref{MainTheorem}, says that we gather all partitions on the first coordinate and the Cantor structure of $C(\lambda)$ to obtain a symbolic model for the attractor $\Gamma_{\lambda}$. This allows us to rigorously prove  several dynamical properties in this case, for instance topological transitivity and the existence of an ergodic invariant measure.

\subsection{Periodic orbits}

The geometric series that arise in the angular dynamics allow us to find periodic orbits with period $m$ and particular itinerary. That is, suppose we are given a word
$w$ of length $m \geq 3$,  given by $w = (w_{j})$ with 
$w_j \in \{+1,-1\} $ for $j=1, \ldots,m$. Each $w_j$
specifies which kind of bounce, left or right, occurs at the $j$th iteration.
If $(s_w,\theta_{w})$ is to be a periodic orbit of period $m$ with this itinerary, then its angle must satisfy the equation
\begin{equation}\label{periodicangle}
\theta_{w} =  \angmap_{w_{m}} \circ \cdots \circ \angmap_{w_1}(\theta_{w}) = 
\frac{\pi}{3} \sum_{j=1}^{m}(-1)^{j-1}w_{m-j+1}\lambda^j + (-\lambda)^m\theta_{w},
\end{equation}
where the second equality gives the explicit formula for the composition.
Solving for $\theta_w$ gives 
$$\theta_{w} = \frac{\pi}{3} \left[ \sum_{j=1}^{m}(-1)^{j-1}w_{n-j+1}\lambda^j \right] \, \left[ \sum_{k=0}^\infty (-\lambda)^{km}\right],$$
since the infinite series equals $[1 - (-\lambda)^m ]^{-1}$.

Setting $b_{km + j} \defeq  (-1)^{km + j-1} w_{m-j+1}$ for any $j=1, \ldots,m$ and $k \in \NN$, we conclude that 
$$\theta_{w} = \sum_{n=1}^\infty b_{n}\lambda^n,$$
% \frac{1}{1-(-\lambda)^m}
and hence $\theta_w \in \elcantor$. 

When $\lambda <\textfrac{1}{3}$, all iterations of the
interval $[0,1]\times \{\theta_{w}\}$ are contained in the Markovian region.
Hence we can partition this interval according to its future itineraries.
Selecting the partition element $J_w$ corresponding to the itinerary $w$, we have $T_\lambda^{m} (J_w) \supset J_w$, so that there is a periodic orbit with periodic itinerary $(w, w, w, \ldots)$. 

\begin{figure*}
\subfigure[Period 4]{%
\includegraphics*[scale=.5]{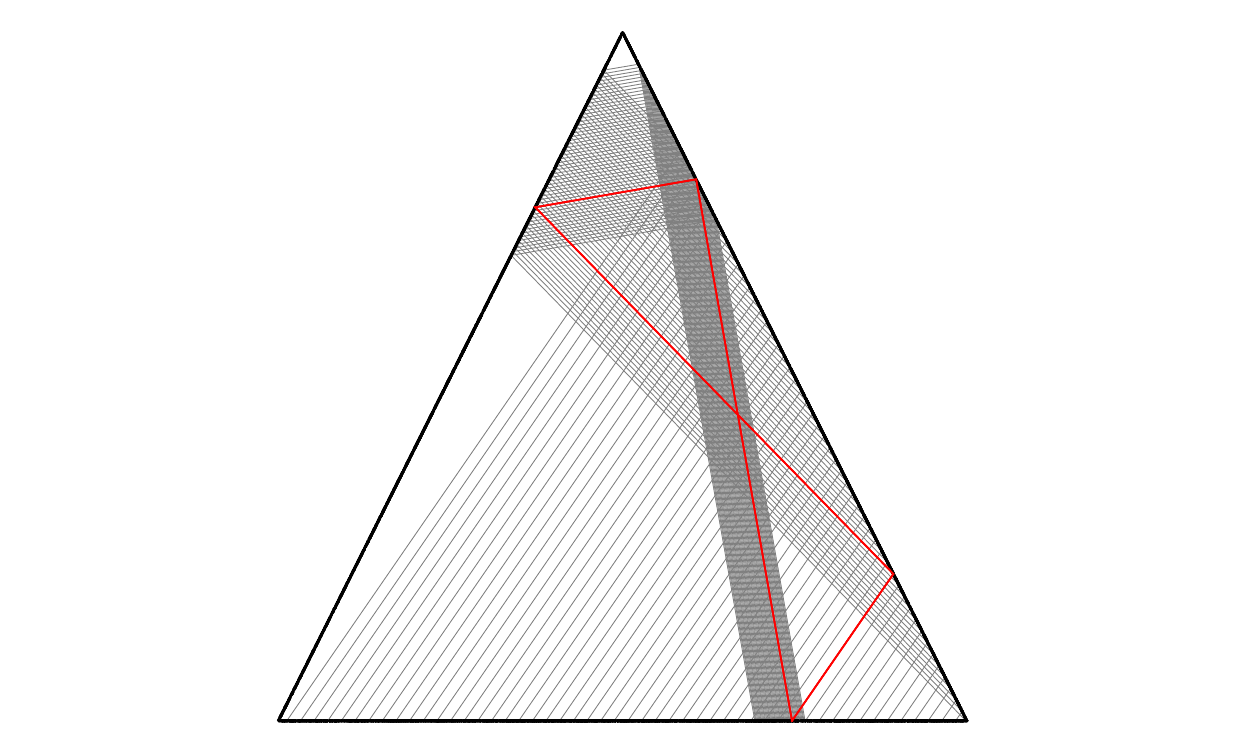} \label{fig:period4}
}
\quad
\subfigure[Period 5]{%
\includegraphics*[scale=.5]{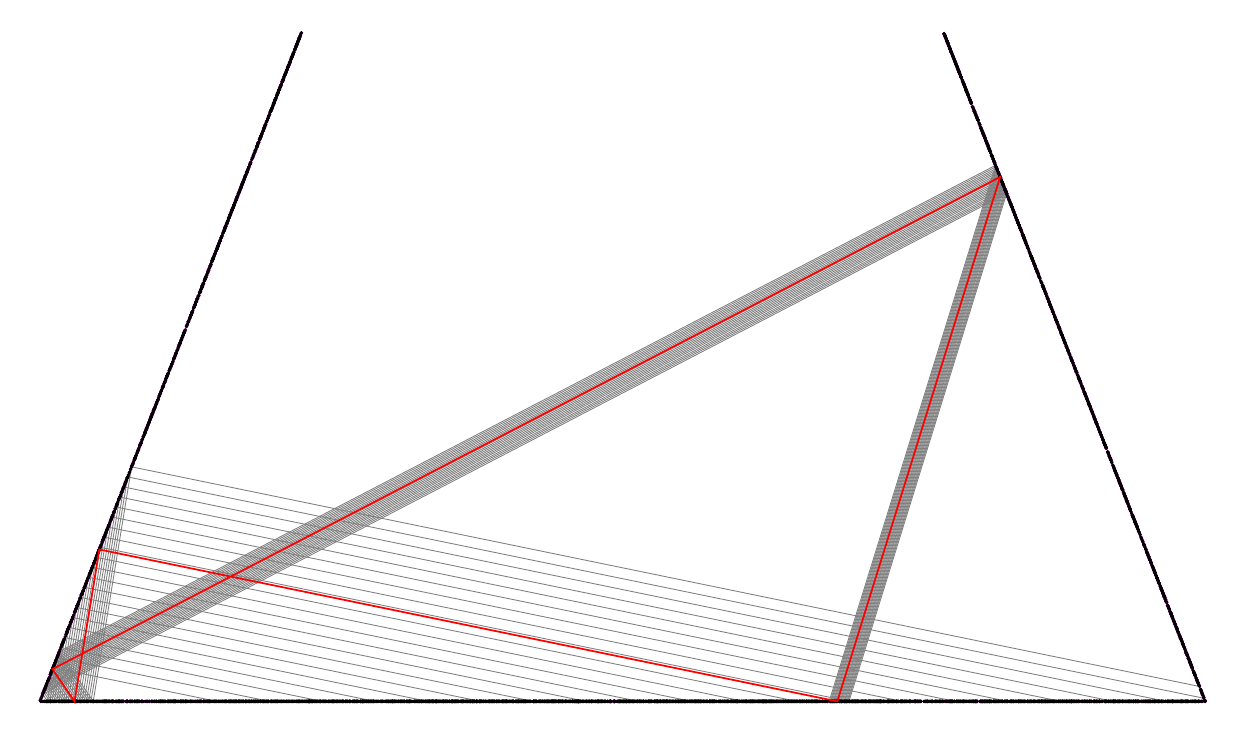}\label{fig:period5}
}
\caption{Periodic orbits (red) in the equilateral triangle for pinbilliard dynamics with $\lambda =0.3$: (a) period-$4$ orbit with itinerary $(+1,+1,-1,-1)$; (b) period-$5$ orbit with itinerary $(+1,+1,+1,-1,+1)$.
The light grey lines show the complete regions with the same itinerary and angles.
 \label{fig:period}}
\end{figure*}

This is the basis of the argument for the rigorous proof of the Main Theorem presented in section \ref{sec:proofMainTheorem}.
Here we exhibit in fig.~\ref{fig:period} two examples of periodic orbits found numerically using this method.
%
%For instance, in  we exhibit two periodic orbits of period $4$ and $5$ for $\lambda= 0.3$, with particular itineraries, and shaded in gray the interval that realizes the corresponding itinerary. 

\section{Non-accessible regions in phase space}

As discussed in the previous section, for $\lambda \le \frac{1}{3}$ the attractor consists of complete horizontal lines, so that for a given angle $\theta \in \elcantor$, (almost) all positions in $[0,3] \times \{\theta \}$ are attained in the attractor.
For $\lambda > \frac{1}{3}$, some of these horizontal lines are cut, giving rise to gaps, located close to the table vertices, that no longer form part of the attractor; see fig.~\ref{fig:phenomenology}(b). These regions are approximately triangular in shape, and change form
as $\lambda$ varies.

To explain this phenomenon, we say that a angle point $(s, \theta)$ is \emph{non-accessible} if $(s, \theta) \notin T_\lambda^n(M),$ for some iterate $n$. 
For example,  
%if $\tilde H = [-\frac{\pi}{2},\frac{\pi}{2}] \setminus [-\lambda \frac{\pi}{2},\lambda \frac{\pi}{2}]$ 
 $(0,3) \times \left(\frac{\lambda \pi}{2}, \frac{\pi}{2}\right) $
 %\cup [0,3] \times \left(-\frac{\pi}{2}, -\frac{\lambda \pi}{2}\right)$$
is such a non-accessible region for $n=1$, along with its  mirror image with negative  angular coordinate. %, denoted by $\LL$. 
Obviously, such non-accessible regions cannot belong to the attractor.

To avoid cumbersome notation, we will mainly work with intervals of angles and their dynamics under $\phi_{\pm 1}$.  We denote by $N \defeq (-\frac{\pi}{2}, \frac{\pi}{2})$ the whole range of available angles,  and  by
$K \defeq (-\frac{\pi}{6}, \frac{\pi}{6})$ the angles corresponding to the Markovian region. Furthermore, given an interval $L = (\ell_1, \ell_2)$ and a scalar $\rho$, we denote  $\rho L \defeq (\rho \ell_1, \rho \ell_2)$. Recall that $\phi_{\pm 1}$ are affine maps of $\RR$ that  depend on 
$\lambda$. 
%Then $\phi_{+1}(N) \cup \phi_{-1}(N) \subset K$.

Simple calculations show that 
the images $\phi_{+1}(K) = (\lambda \frac{\pi}{6}, \lambda \frac{\pi}{2})$ and $\phi_{-1}(K) = -\phi_{+1}(K)$ are disjoint intervals, both contained in 
$\lambda N$, and the gap between them is $\lambda K$; that is, we have
$\lambda N = \phi_{+1}(K) \cup \lambda K \cup \phi_{-1}(K)$,  together with two points at the boundaries of the intervals.

If $\lambda < \frac{1}{3}$, then we have that $\lambda N \subset K$, and hence we can conclude that
the strip $(0,3) \times \lambda K \subset M$ is a non-accessible region, since it is not contained in $T_\lambda(M)$.
Iterating this procedure gives rise to the Cantor set $\elcantor$ discussed in the previous section.
%; note that the gap in the $n$th step of the construction is of order $\lambda^{n}$.

When $\lambda>\frac{1}{3}$, we now have the opposite case:
$\lambda N \supset K$, so there are points that may enter the middle gap. Since we have studied the images of points inside $K$ under the action of $\phi_{+1}$ and $\phi_{-1}$, it remains to study the images of 
$\lambda N \setminus K$.

Let $\alpha \defeq \phi_{+1}(+\lambda \frac{\pi}{2}) = \textfrac{\pi}{6}(2 \lambda - 3\lambda^2) $; by symmetry, we have $\phi_{-1}(-\lambda \frac{\pi}{2}) = -\alpha$.
For $ \lambda > \frac{1}{3}$, we have $\alpha < \lambda \frac{\pi}{6}$, so that the image of $\lambda N \setminus K$ covers part of the previously-found gap. 
When $\alpha = 0$, at $\lambda = \frac{2}{3}$, the gap is completely covered by this image.

The images of the interval $\lambda N \setminus K$ consist of two symmetric intervals, $L \defeq (\alpha, \lambda \frac{\pi}{6})$ and $-L$.  
It is, however, necessary to study in more detail the position coordinate of the image under $T_\lambda$ of the strip $(0,3) \times L$, since $\lambda N \setminus K$ is no longer in the Markovian region, and so the image of horizontal intervals of the form  $(i, i+1) \times \theta$ do not completely cover a side.

We thus only need to study the image of the rectangle $R \defeq [0,1]\times (\frac{\pi}{6}, \lambda \frac{\pi}{2}) \subset M^+_{0}$, where $M^+_0$ is in the domain of $\phi_{+1}$. 
We have already shown that the image $T_\lambda(R) \subset [1,2] \times L$.

In order to obtain a better notion of the shape of the image, we have to be careful about the images of each of the four sides of the rectangle $R$.  The map $T_{\lambda}$ is, of course, not defined on vertical sides, but we  can take its  continuous extension.

The bottom side of $R$, given by $(0,1)\times \{ \frac{\pi}{6} \}$, is mapped onto the line $(1,2)\times \{ \lambda \frac{\pi}{2} \}$, while the top side, 
%On the other hand, the upper horizontal side of the rectangle, \
$(0,1)\times \{\lambda \frac{\pi}{2}\}$, is mapped onto a sub-interval $(1,\ell)\times \{  \alpha \}$, reversing the orientation. The point $\ell = \ell(\lambda)$ corresponds to the point where a ball hits the triangle's boundary if it starts on the vertex $0$ with angle $\lambda \frac{\pi}{2}$. 

For the right side of $R$, it is necessary to consider vertical lines $\{1-\epsilon\} \times (\frac{\pi}{6}, \lambda \frac{\pi}{2})$ for small $\epsilon>0$, since the image is not actually defined for the right side of $R$.
 It is clear that the projection on the angular coordinate of the image covers the whole interval $\left( \alpha ,\lambda \frac{\pi}{6} \right)$, since the angles do not depend on the position $s$ on the table. Furthermore, as $\epsilon \to 0$, both ends of the curve converge to the vertex at $s=1$. So the continuous extension of $T_{\lambda}$ maps the right side of the rectangle onto the vertical interval $\{1\} \times \left(\alpha ,\lambda \frac{\pi}{6} \right)$. 

Finally, the left side of $R$, given by $\{0\} \times (\frac{\pi}{6}, \lambda \frac{\pi}{2})$, is mapped to a curve joining the point $(2,\lambda \frac{\pi}{2})$ and the point $(\ell, \alpha)$.
In fact, this curve is given by $\gamma(t) \defeq (\ell(t), \phi_{+1}(t))$ for $t\in (\frac{\pi}{6},\lambda \frac{\pi}{2})$, where 
\begin{equation} \label{eq:cuttingline}
\ell(t) \defeq 1 + \frac{2}{1 + \sqrt{3} \tan(t)}.
\end{equation}
Note that only the endpoint  $\ell = \ell(\lambda \frac{\pi}{2})$, and not the curve itself, depends on $\lambda$.
%  \aubin{Es posible que no sea relevante escribir esta f\'ormula. Adem\'as, la f\'ormula que estaba en la version 30jul corresponde a otra linea de corte.}

%The complement in $M_{0}^+$ of the rectangle $R$ studied  above is mapped into the set $[1,2]\times(  \lambda \frac{\pi}{6},  \lambda \frac{\pi}{2})$.  Therefore, 
To summarise, we have shown that $T_{\lambda}(M_{0}^+)$ is contained in $M_{1}$ in the region bounded below by the curves $\max\{\alpha, \gamma(t)\}$, and above by the line $\lambda \frac{\pi}{2}$.

\begin{figure*}
\subfigure[$\lambda=0.3$]{%
\includegraphics*[scale=.15]{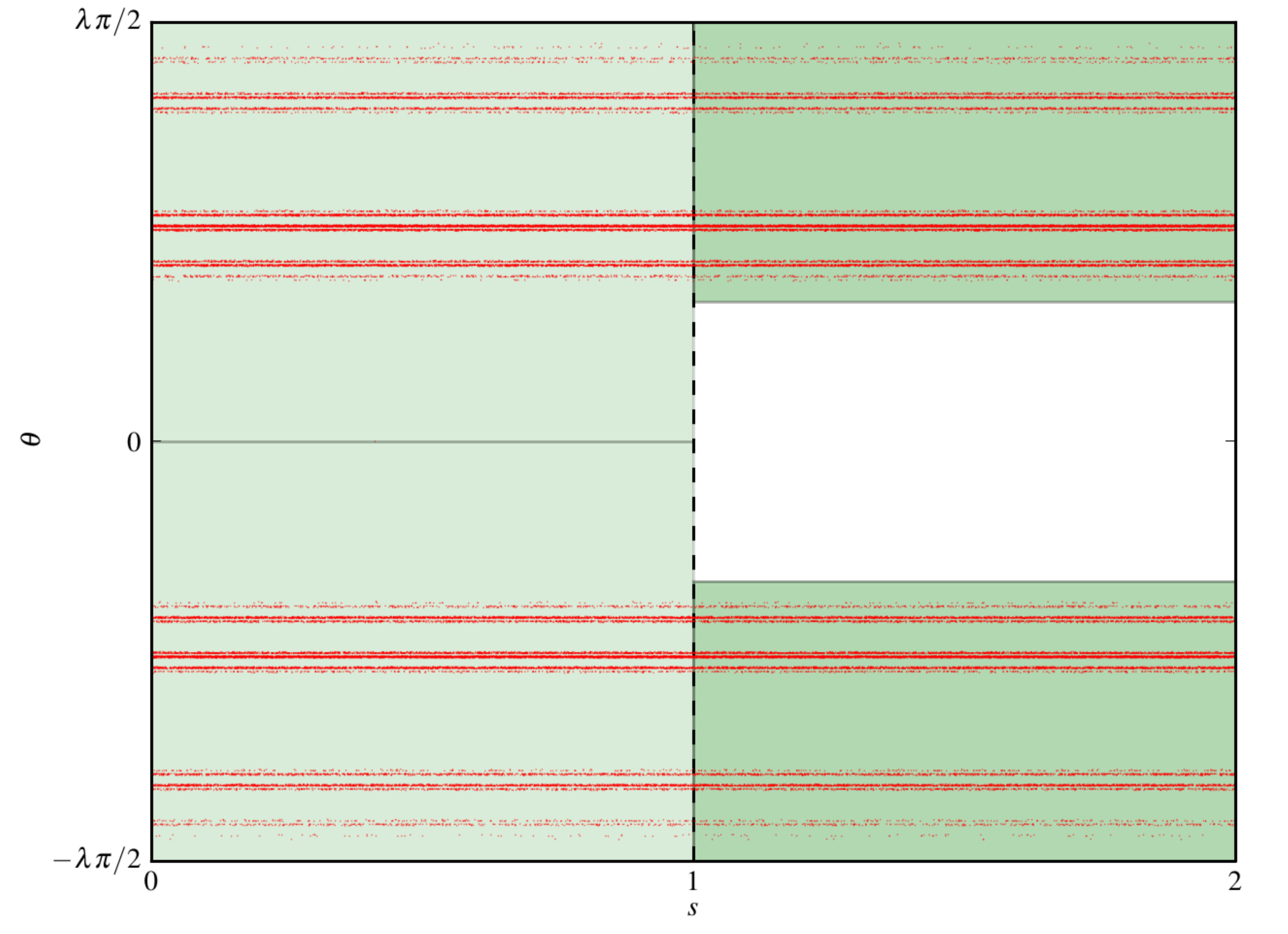}
}
\subfigure[$\lambda=0.5$]{%
\includegraphics*[scale=.15]{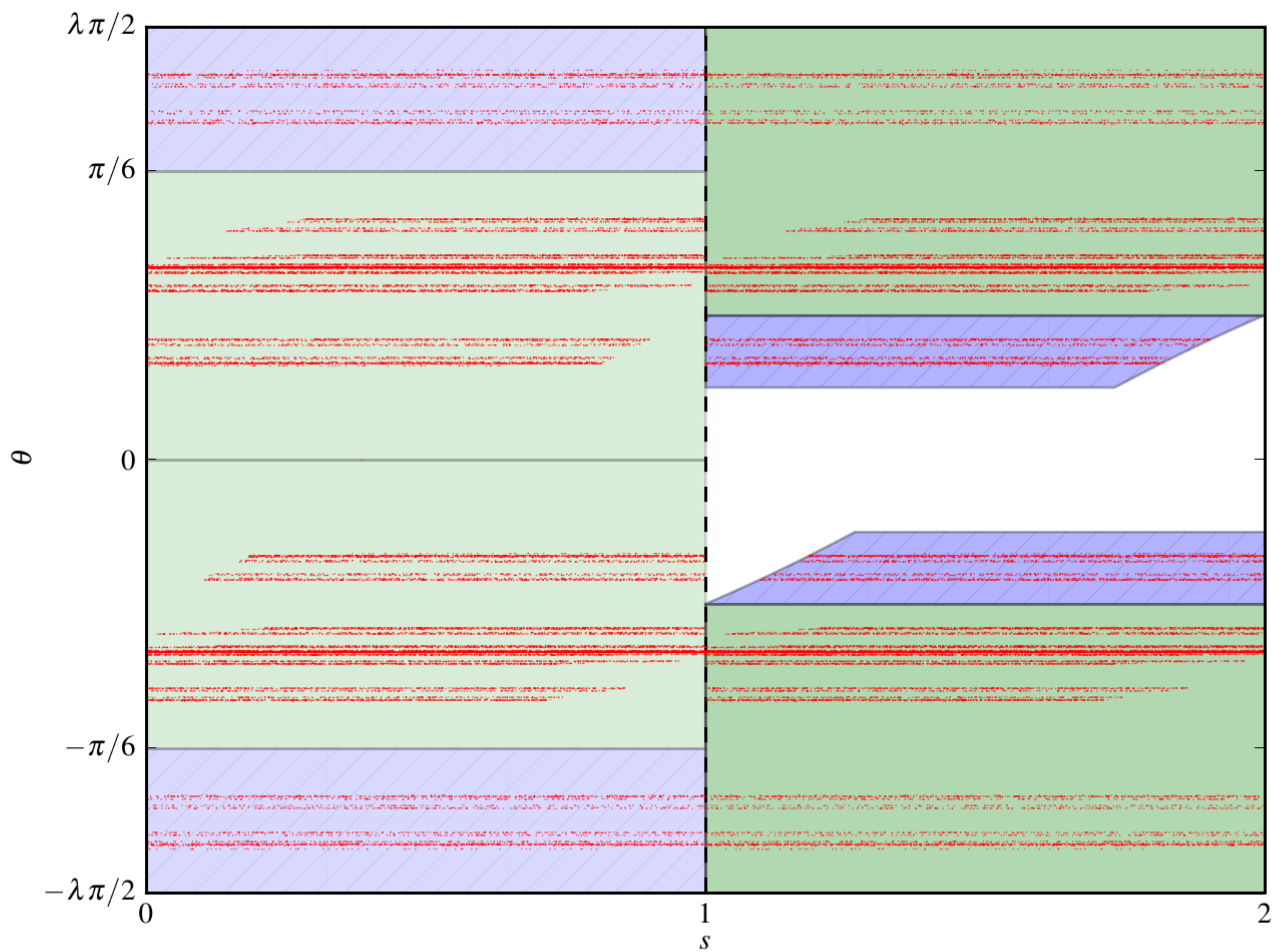}
}
\subfigure[$\lambda=0.75$]{%
\includegraphics*[scale=.15]{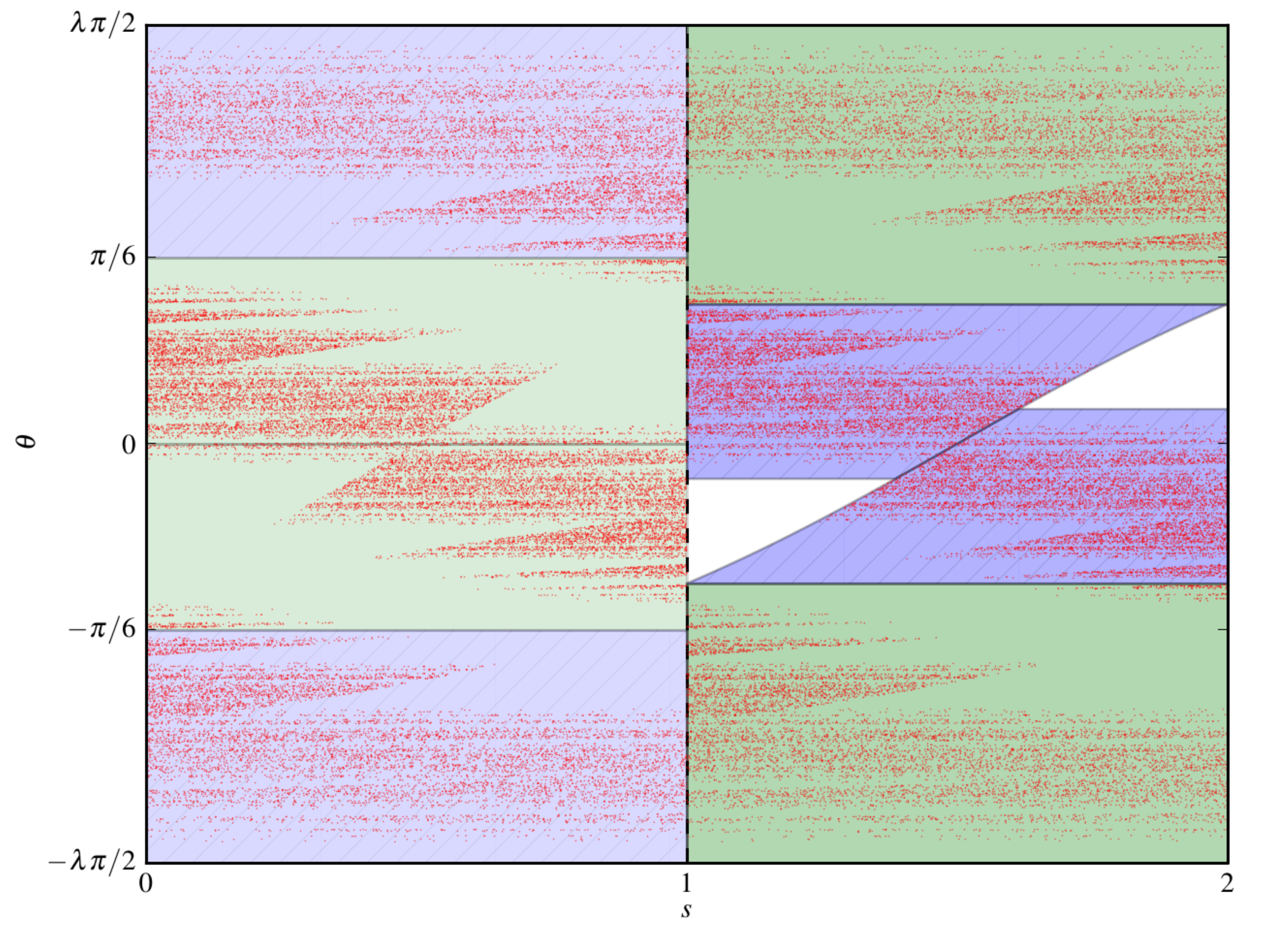}
}
\caption{Action of $T_\lambda$ on $\lambda M_0$ for different values of $\lambda$.
The Markovian region is shaded in green and the non-Markovian region in blue. 
Darker shading is used for the images, shown for $s \in (1,2)$; lighter shading is used for the domain, shown for $s \in (0,1)$. The attractor is shown in red dots. }
\label{fig:lineasdecorte}
\end{figure*}

The other way that trajectories can be injected into $M_{1}$ is if they come from $M_{2}^-$. By symmetry, the relevant image of $M_{2}^-$ is bounded by the curves 
$-\lambda \frac{\pi}{2}$ and $\min\{-\alpha, \gamma(t)\}$. This is shown in figure \ref{fig:lineasdecorte}, and explains the presence of the triangular-shaped gaps in the attractor.
\

As we have seen, when $\lambda > \frac{2}{3}$, 
the angular parts of the images overlap. In fig.~\ref{fig:lineasdecorte}(c), we see that the complete images in fact touch. Although this makes possible that the central bands of the attractor  join  as soon as $\lambda > \frac{2}{3}$, in fact we observe numerically that this occurs later, around $\lambda \simeq 0.7$. This can be understood by tracing the evolution of the non-accessible regions in time and finding the limits of the region around $\theta=0$ which is correspondingly excluded from the attractor.

\section{Isolation of period-3 orbit from the attractor}
Recall from figure \ref{fig:phenomenology}(b) that a splitting of the central bands of the attractor is observed around $\lambda \simeq 0.68$. This is related to the  
 local stable and unstable manifolds of the distinguished period-$3$ orbit discussed in section~\ref{3per}.

\begin{figure*}
% \subfigure[$\lambda=0.3$]{%
% \includegraphics*[scale=.5]{equilateral_lambda0_3} \label{fig:triangle:a}
% }
% \quad
\subfigure[$\lambda=0.5$]{%
\includegraphics*[scale=.45]{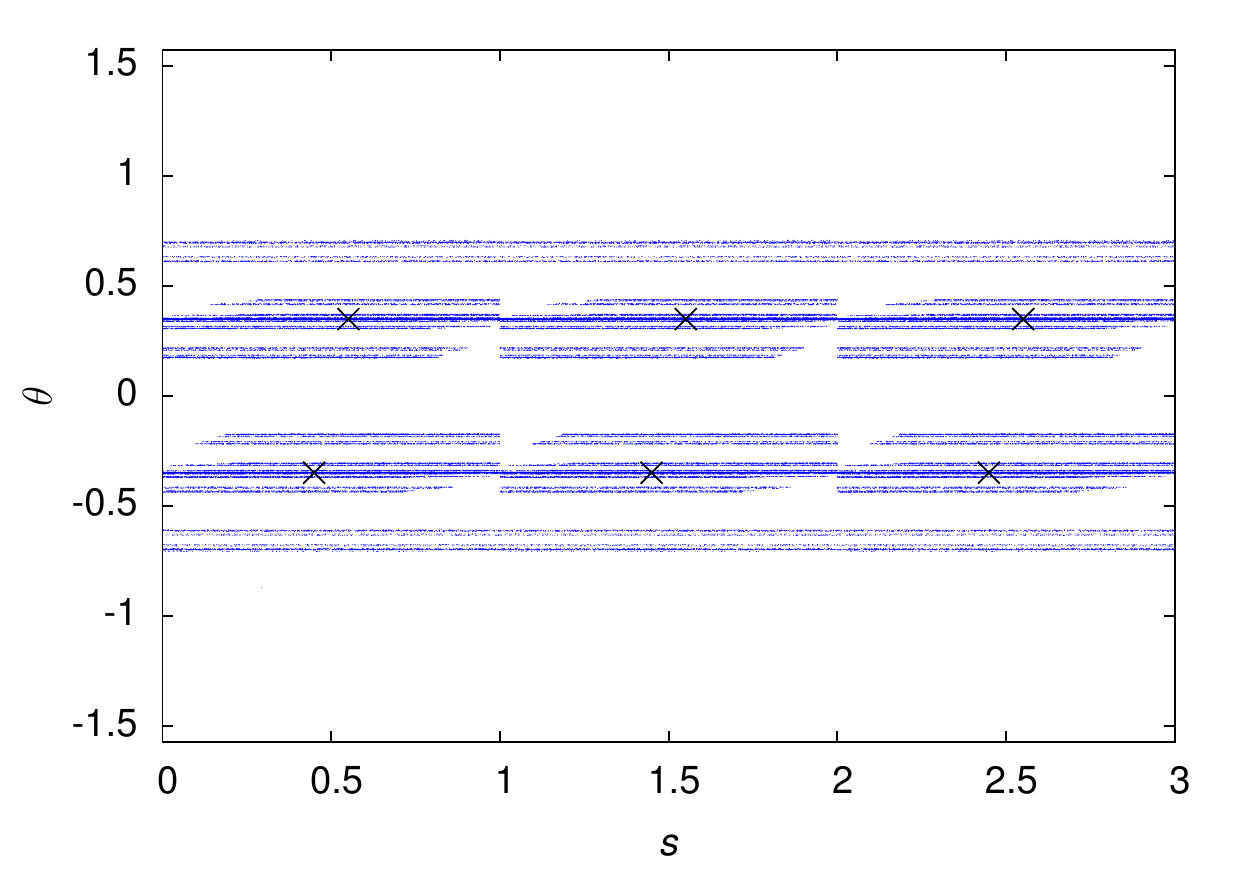}
}
\subfigure[$\lambda=0.68$]{%
\includegraphics*[scale=.45]{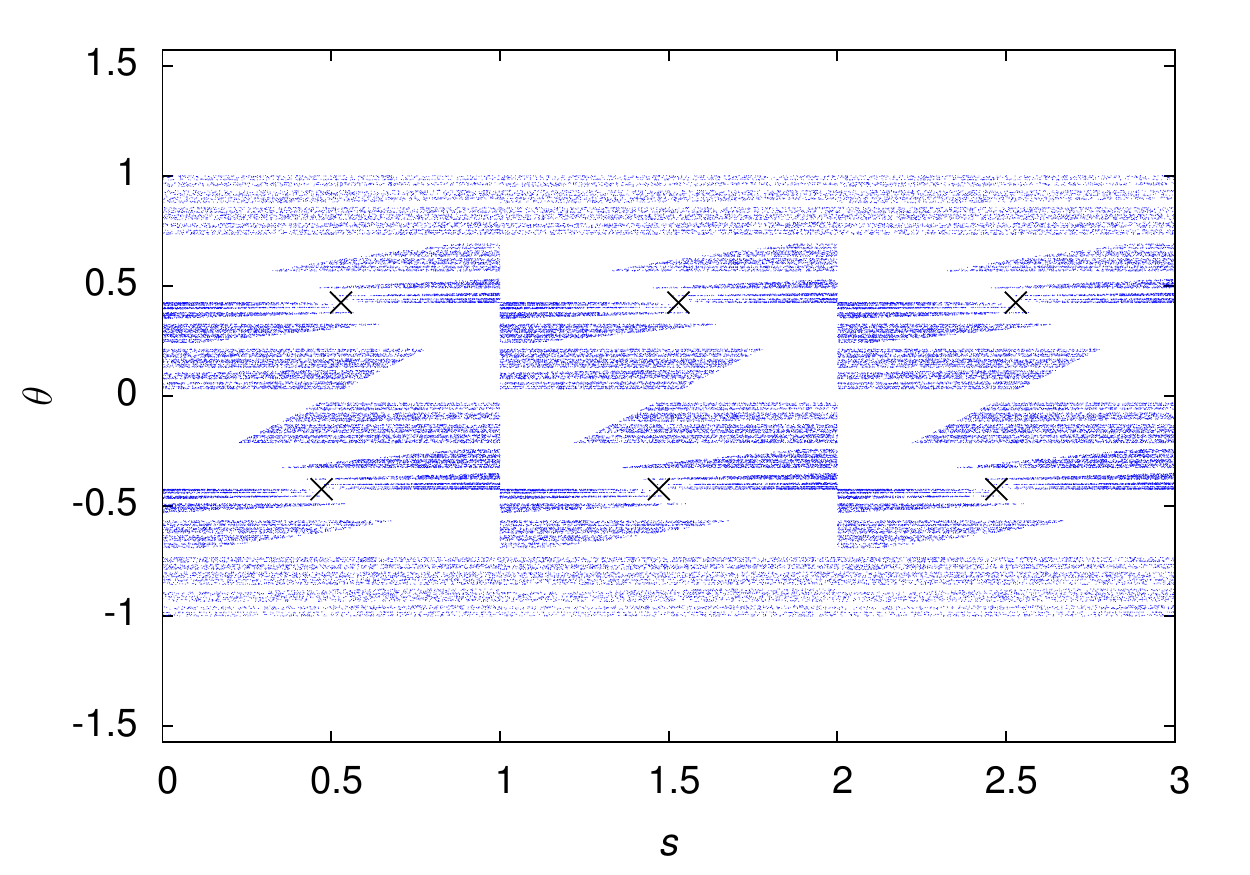}
}
% \quad
\subfigure[$\lambda=0.8$]{%
\includegraphics*[scale=.45]{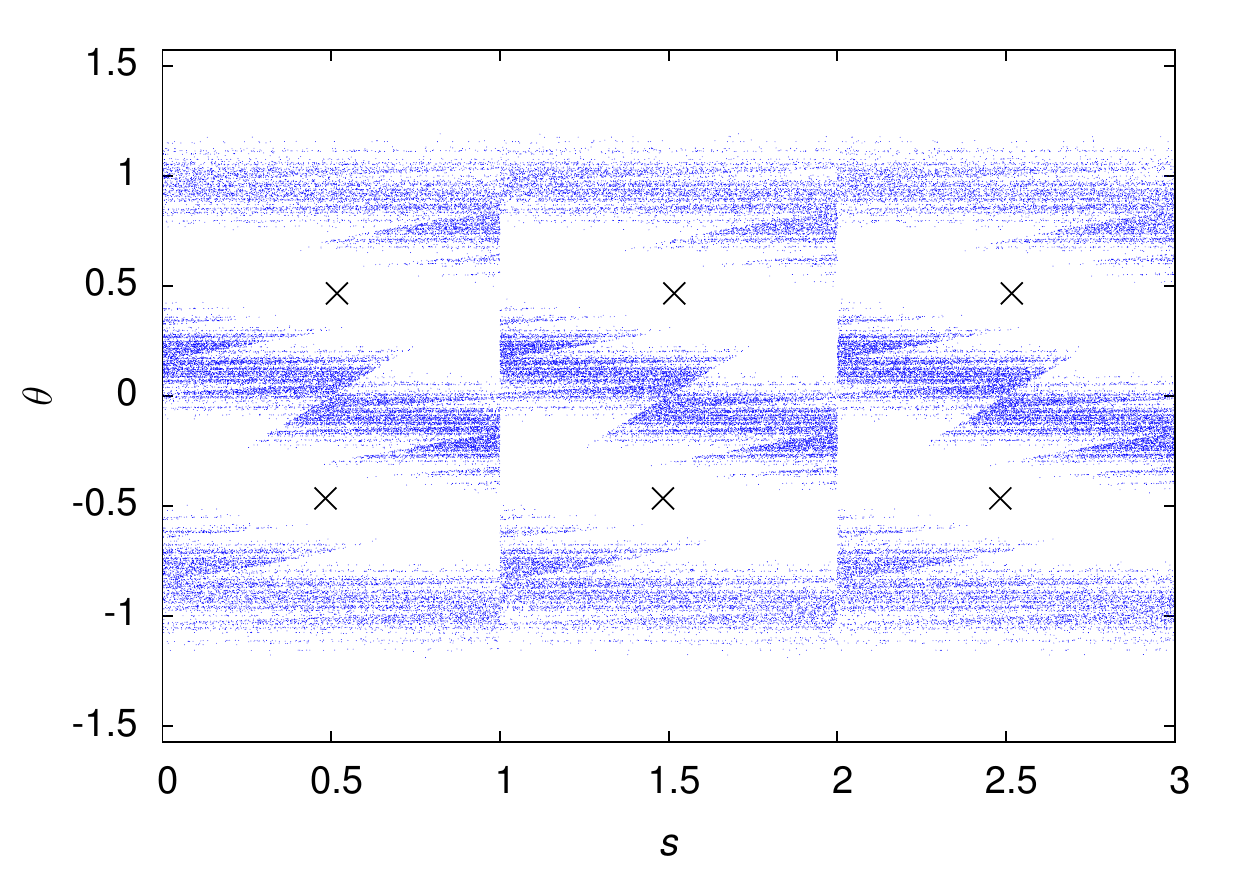}
}

% \qquad
% \subfigure[$\lambda=0.95$]{%
% \includegraphics*[scale=.5]{equilateral_lambda0_95}
% }
% \quad
% \subfigure[$\lambda=0.98$]{%
% \includegraphics*[scale=.5]{equilateral_lambda0_98_espacio_fase_orig}
% }
% 
\caption{Phase space of pinbilliard dynamics in the equilateral triangle, with the period-$3$ orbits shown as crosses.
% Plots of $10^5$ iterates of a single point for  $\lambda =0.
}
\label{fig:period-3-bifurcation}
\end{figure*}

Figure~\ref{fig:period-3-bifurcation} shows the period-$3$ orbit together with the attractor. From the figure, we see that  the period-$3$ orbit is embedded in the attractor before the bands split, and in particular that its local unstable manifold
is a complete horizontal line contained in the center of the band.  After the band splits, however, the period-$3$ orbit becomes isolated from the strange attractor, and there is no longer such a central line.

%\subsection{When period 3 orbits leave the attractor}
% There is one bifurcation phenomenon that occurs for some parameter $0.55< \lambda_{b} <0.75$, and this is when the period 3 orbits get isolated from the rest of the attractor. 

\begin{figure*} %  figure placement: here, top, bottom, or page
\subfigure[$\lambda=0.55$]{
    \includegraphics*[scale=0.6]{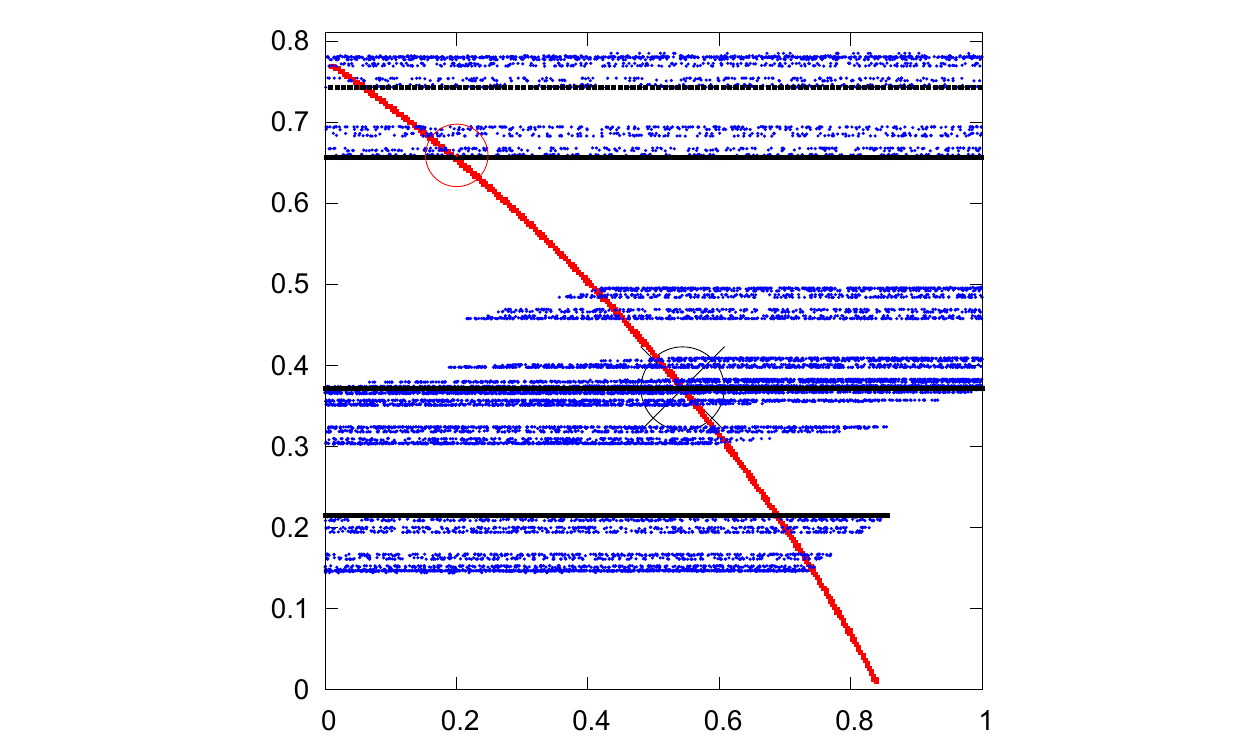}
}
\subfigure[$\lambda=0.75$]{
    \includegraphics*[scale=0.6]{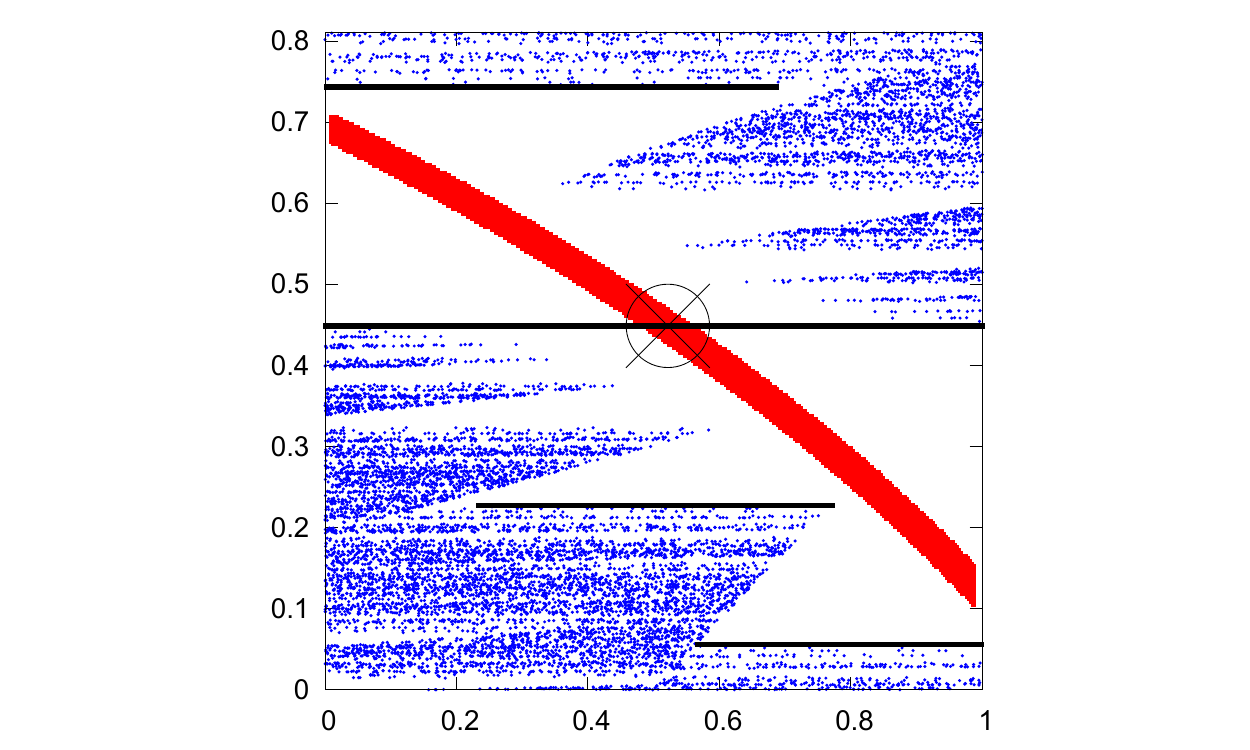}
}
\caption{The region $M_0$ with non-negative angles, showing the strange attractors (blue points), one point in the period-$3$ orbit (crossed circle), several iterations of the unstable manifold of the periodic point (thick black line), and a region in which its local stable manifold must be contained (red area) for (a) $\lambda = 0.55$ and (b) $\lambda=0.75$. In (a), there is a homoclinic point (intersection of the stable and unstable manifolds), marked with a small red circle. A number of iterations $n=30$ was used. }
    \label{fig:homoclinicpoint}
 \end{figure*}

This phenomenon is due to a homoclinic bifurcation: for $\lambda = 0.55$, the stable and unstable manifolds of the period-$3$ meet, forming a homoclinic point, as shown in fig.~\ref{fig:homoclinicpoint}(a). 
Thus any neighbourhood of the period-$3$ orbit contains recurrent points, and hence the period-$3$ point is not isolated from the dynamics.

By $\lambda \simeq 0.75$ the stable and unstable manifolds no longer meet, and this recurrence is destroyed, allowing the period-$3$ orbit to become isolated.  

To obtain these figures, the unstable manifold is obtained by iterating a local horizontal interval around the period-$3$ point. 
The stable manifold, on the other hand,  is approximated using the following escape-time algorithm. For a given point in phase space to belong to the local stable manifold of the period-$3$ orbit, its future iterates must follow the same dynamics as the period-$3$ orbit, e.g.\ only have right bounces. 
We thus iterate   each point in a grid  on the domain a certain number $n$ of times; if the point has only right bounces up to this time, then the initial condition is marked as forming part of the approximate local stable manifold of order $n$. Of course, this does not give the precise shape of the stable manifold, but rather a region in which it must be contained. In particular, this guarantees that there is no homoclinic point if the unstable manifold does not cross this region.

This phenomenon also provides evidence of the relation  \cite[Section~III.3]{CM03} between
expansion in the unstable direction and fractioning due to the set $\DiscontinuityRegion$. 
The expansion in  the unstable direction of the  
period-$3$ orbit is given by
$[\frac{\cos (\lambda \theta_*)}{\cos(\theta_*)}]^3$, which  decreases as
$\lambda$ increases.
%; for a general discussion, see, for example,
%.
  
%\section{Band joining}
%Although the images of phase space touch and join each other for $\lambda=\frac{2}{3}$, the bands in the attractor remain separate until around $\lambda = 0.7$. Apparently the attractor does not actually attain the limits of the region?

\section{Non-transitive behaviour close to $\lambda = 1$} \label{sec:nontransitive}

\begin{figure*}[tp] %  figure placement: here, top, bottom, or page
   \centering
   \subfigure[]{
   \includegraphics[scale=0.5]{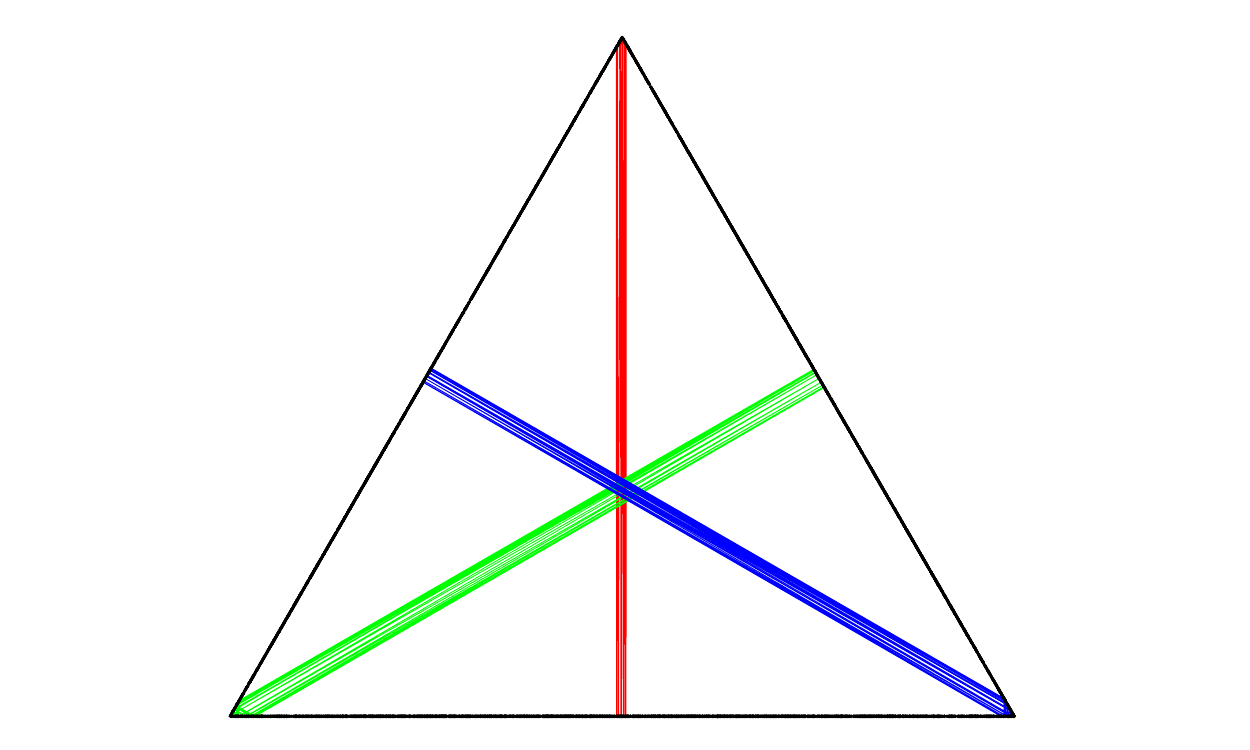}
   }
   \quad
   \subfigure[]{
   \includegraphics[scale=0.2]{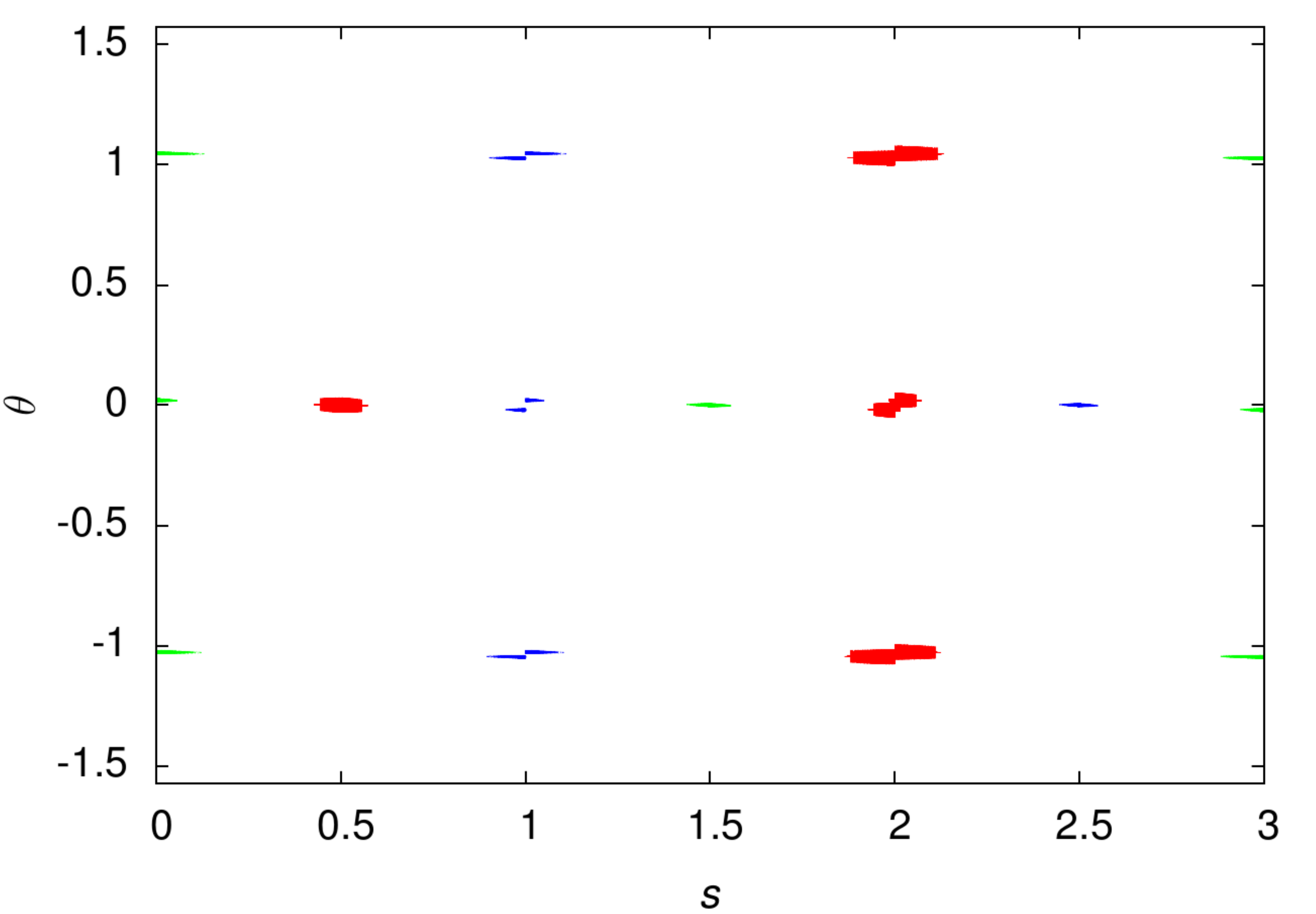}
   }
   \caption{
Three non-communicating transitive parts of the attractor, shown in different colors for $\lambda = 0.98$, in (a) configuration space and (b) phase space.   
   \label{fig:nontransitive}
   }
\end{figure*}

For $\lambda \simeq 0.98$, we observe that there are three distinct transitive parts of the attractor, i.e.\ such that it is impossible to jump from one to the other; see fig.~\ref{fig:nontransitive}.

The  reason for this phenomenon is the following.
After enough iterations, the system arrives close to a vertex, with an angle close to $\frac{\pi}{3}$, so that it collides with the opposite side at an almost perpendicular angle, and is reinjected back into the same corner.
This acts as a trapping region which prevents the system from escaping to one of the other parts of the attractor.

\begin{figure*}[tp] %  figure placement: here, top, bottom, or page
   \centering
   \subfigure[]{
   \includegraphics[scale=0.15]{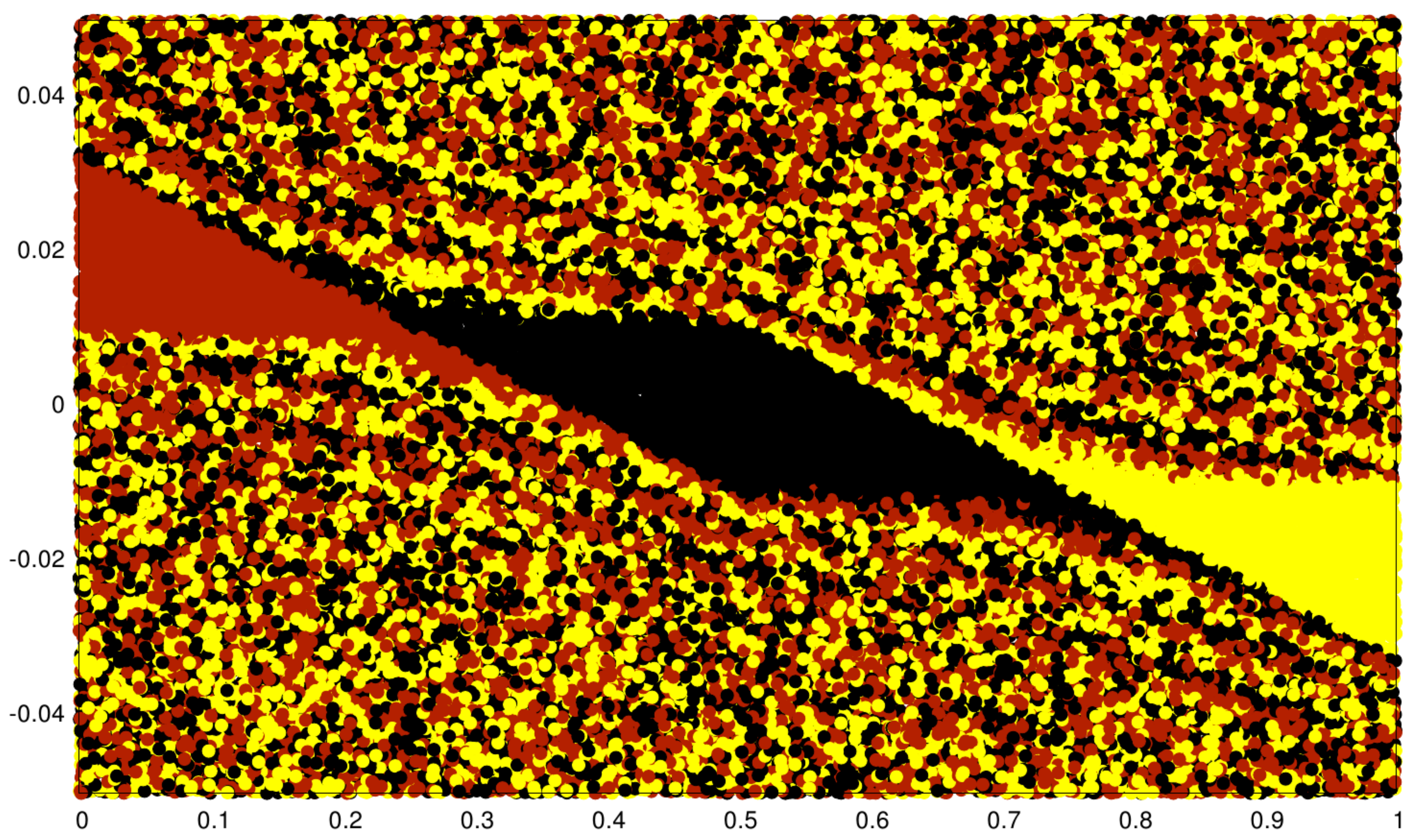}
   }
   \quad
   \subfigure[]{
    \includegraphics[scale=0.15]{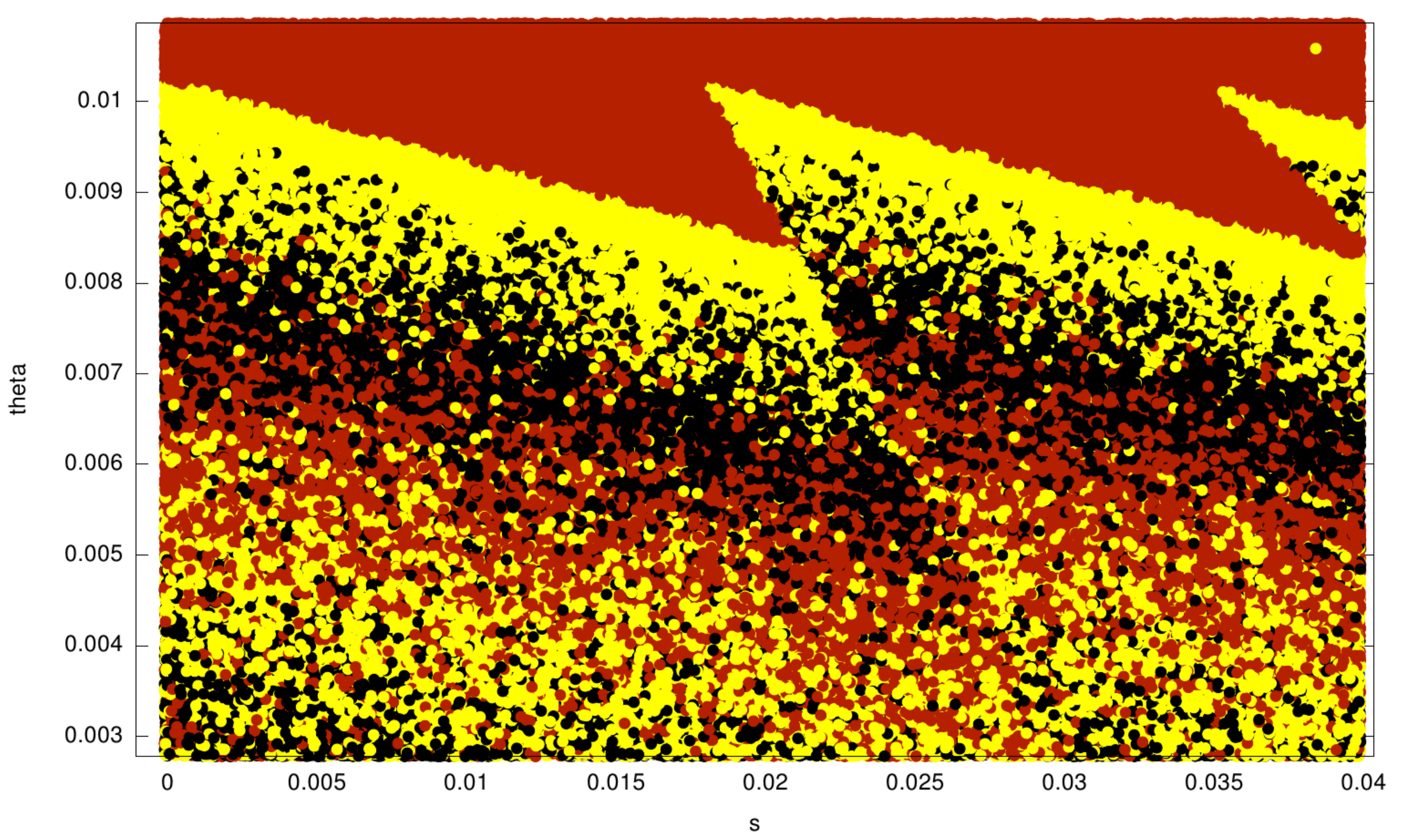}
   }
   
   \caption{%
Basins of attraction of the three transitive parts of the attractor, shown in three different colors, exhibiting a fractal basin boundary.  (a) Region around one part of one of the transitive parts of the attractor, centered at $(0.5, 0)$. (b) Detail of the basin boundary.  }
   \label{fig:nontransitivebasins}
\end{figure*}

In~\ref{fig:nontransitivebasins} the basins of attraction of each of the three symmetrical copies of this region are shown, indicating a fractal basin boundary structure.

%%%%%%%%%%%%%%%%%%%%%%%%%

\section{Proof of Main Theorem} \label{sec:proofMainTheorem}
In this section we give a rigorous proof of the Main Theorem. This will be done in two parts, treating separately the topological and the ergodic properties, although both are consequences of the existence of a topological conjugacy of the pinball billiard map in the attractor with a symbolic dynamics model. 

Let us recall that the attractor of $T_{\lambda}$ (or maximal forward invariant set in $M$) is defined by 
$\maximalinvariantset \defeq \cap_{n\geq0}T_{\lambda}(M) \subset M$. 
Roughly speaking, the symbolic model will be a shift-like map in the space of bi-infinite sequences of $\pm1$, but interchanging signs in a certain part of the sequence at each iteration; this map is defined in eq.~\eqref{shiftmodel}. This map possesses a natural ergodic measure which will be pulled down to an invariant, ergodic measure for the pinball billiard map using the above-mentioned conjugation. 

The Markovian property that guarantees uniform expansion along certain horizontal intervals is the cornerstone in the construction of the conjugation, discussed in sec.~\ref{sec:markovian-heuristic} and constructed in section \ref{markovianReg} below. In some cases, for instance for $0<\lambda < \frac{1}{3}$, the attractor is contained in the region where this property holds, allowing us to prove the existence of a dense orbit.
 
It is interesting that the set of geometric series with random signs  is intimately related with the action of the collision rule of the pinball billiard map on the angular coordinate. 

Let  $0<\lambda <1$ and recall the set $\elcantor$ defined in \eqref{ElCantor}:
$$\elcantor \defeq \left\{ \frac{\pi}{3} \sum_{n=1}^\infty z_n \lambda^n :  z_n \in \{1,-1\} \, \forall n \geq 1 \right\}.$$ 
%%%%%%%%%%%%%%%%%%%%%%%%%%%%%%%%%%%%%%%%

%$$\elcantor = \left\{ \frac{\pi}{3} \sum_{n=1}^\infty w_n \lambda^n | \text{ where } w_n \in \{1,-1\}, \, \forall n \geq 1 \right\} \subset \frac{\pi}{3}\left[ -\frac{\lambda}{1-\lambda}, \frac{\lambda}{1-\lambda} \right]. $$

This is a subset of the interval $\frac{\pi}{3}[-\frac{\lambda}{1-\lambda},\frac{\lambda}{1-\lambda}]$. It is a Cantor set if $\lambda < \textfrac{1}{2}$, and all numbers in $C(\lambda)$ then have a unique expression in terms of a signed geometric series, whereas for $\lambda \geq \textfrac{1}{2}$ it is a full interval, which grows to fill up $\RR$ as $\lambda \to 1$.  

To prove the Cantor structure of this set, note that if $\lambda < \textfrac{1}{2}$ then for each integer $m \geq 1$ we have
$$  \delta_{m}:=\frac{\pi}{3}\lambda^m - \frac{\pi}{3}\sum_{j=m+1}^\infty \lambda^j  = \frac{\pi}{3}  \, \frac{ 1-2 \lambda }{ 1-\lambda }  \lambda^m   > 0.  $$
Hence, there are no points of $\elcantor$ in any open interval of the form $(s-\delta_{m}, s+\delta_{m})$, for any partial sum $s = \frac{\pi}{3} \sum_{j=1}^{m-1} \pm \lambda^j$. In the same fashion, the length of each bounded connected component of the complement of the union of these intervals is of order $\lambda^{m+2}$.

On the other hand, if $\lambda < \textfrac{5}{6}$ then $\elcantor \subset [-\textfrac{\pi}{2}, \textfrac{\pi}{2}]$. Recall that the set $C(\lambda)$ is invariant under both maps $\angmap_{\pm 1}$, i.e.\ $\angmap_{\pm1}(\elcantor) \subset \elcantor$. 

To observe the relation 
between the geometric series and the pinbilliard,
consider the topological space
$$\elshift \defeq \{ (z_n) : z_n \in \{1,-1\} \, \forall n \ge 1 \}, $$ with the standard product topology,
and the modified shift transformation $-\sigma : \elshift \to \elshift$ given by  $-\sigma((z_{n})) := (- z_{n+1})$. In these terms,  if $0<\lambda<\textfrac{1}{2}$, then both maps $\angmap_{\pm1}$  are branches of the inverse $ E_{\lambda}^{-1}(-\sigma(z))$,
%, for any $z \in \elshift$, 
where the map $E_\lambda: \elshift \to \elcantor$ is defined by
\begin{equation} \label{erdosmap}
E_\lambda(z) \defeq \frac{\pi}{3}\sum_{n= 1}^\infty z_n \lambda^n.
\end{equation}
The condition $0<\lambda <\textfrac{1}{2}$ guarantees that $E_{\lambda}^{-1}: \elcantor \to \elshift$ is well-defined.

On the other hand, to any point in phase space having an entire well-defined future orbit there corresponds a sequence $(w_n) \in \elshift$, such that $w_n = +1$ if the $n$th iterate is a right bounce, or $w_n = -1$ if it is a left bounce. So, if $(s_0,\ang_0) \in M \setminus \DiscontinuityRegion$, then for any $n\geq 1$ we have $T^n(s_0,\ang_0) = (s_n,\ang_n)$, where
\begin{equation}\label{AngleFormula}
\ang_n =  \angmap_{w_{n}} \circ \cdots \circ \angmap_{w_1}(\ang_0) = 
\frac{\pi}{3} \sum_{j=1}^{n}(-1)^{j-1}w_{n-j+1}\lambda^j + (-\lambda)^n\ang_0.
\end{equation}

In this way, we can define the function $\itin : \maximalinvariantset \to \elshift$
by $\itin(P) := (w_n) \in \elshift$ for $P\in \maximalinvariantset$.
If a point hits a vertex of the table at some iterate $k\geq 1$ then one cannot  compute its itinerary after the $k$th step. In this case, we take the convention that the tail of the itinerary is either $w_{k+j}=(-1)^j$ or $(-1)^{j+1}$, for $j \geq 1$.
% Therefore, the map $\itin: M \to \elshift$ is well defined for any point in phase-space.

Recall that for any $0<\lambda<1$ we have $T_{\lambda}(M) \subset (0,3) \times 
(-\lambda \textfrac{\pi}{2}, \lambda \textfrac{\pi}{2})$, and note that if $0<\lambda<\textfrac{1}{3}$ then $\elcantor \subset (-\lambda \textfrac{\pi}{2}, \lambda \textfrac{\pi}{2})$. For such $\lambda$, the attractor $\limitset$ inherits the geometry of $\elcantor$, as shown in the following lemma. 

\begin{lemma} \label{iscontained}
If $0<\lambda<\textfrac{1}{3}$ then $\limitset \subset (0,3) \times \elcantor$.
\end{lemma}

\proof Consider any point $(q,\ang) \in \limitset $. Then there exists a point $(p_0,\ang_0) \in M \setminus \DiscontinuityRegion$, and an increasing sequence of integers $n_k \to \infty$, such that $ (p_{n_k},\ang_{n_k}):= T_{\lambda}^{n_k}(p_0, \ang_0)  \to (q,\ang)$ as $k\to \infty$. Using the expression in \eqref{AngleFormula}, we have, for any $k \geq 1$, that
\begin{eqnarray*}
\left| \ang_{n_k} - \frac{\pi}{3} \left( \sum_{j=1}^{ n_k }(-1)^{j-1}w_{n_k-j+1}\lambda^j  +
  \sum_{j=n_k + 1}^{ \infty } \pm \lambda^j \right) \right| \\
<  \lambda ^{n_k}\frac{\pi}{2} + \frac{\pi}{3}\frac{\lambda^{n_k+1}}{1-\lambda} < C \lambda^{n_k},
\end{eqnarray*}
for some positive constant $C$
%=\pi(1/2 + \lambda/(3-3\lambda))$;
and any selection of $\pm1$ in the series on the right-hand side. Thus $\ang \in \elcantor$. \qed

\subsection{Region with  a Markovian property} \label{markovianReg}
The discontinuity curve $\delta$, defined in \eqref{eq:curvasdiscontinuidad}, does not depend on the parameter $\lambda$, and its image is contained in the rectangle $(0,3)\times(-\textfrac{\pi}{6}, \textfrac{\pi}{6})$. Furthermore, phase space is naturally partitioned into horizontal intervals of the form $\hinterval{i}{\ang} := (i,i+1) \times \{\ang \}$ for $i \in \ZZ_3$ and $\theta \in (-\textfrac{\pi}{2},\textfrac{\pi}{2})$. We denote this partition by $\mathcal{J}$. 

We define $\MarkovianRegion \subset \mathcal{J}$ as $\MarkovianRegion \defeq
\{ J(i, \theta) \in \mathcal{J} : |\theta| < \frac{\pi}{6} \}$.
The region of phase space covered by $\MarkovianRegion$, i.e.\ the union of all $\hinterval{i}{\ang} \in \mathcal{J}$ such that $|\ang|<\textfrac{\pi}{6}$, has a non-uniform expanding behavior along horizontal lines.
%The \emph{markovian} region of the pinbilliard correspond to the union of all $\hinterval{i}{\ang} \in \mathcal{J}$ that $|\ang|<\pi/6$; let us denote this subset of $\mathcal{J}$ by $\MarkovianRegion$. 
To see this, note that any $J \in \MarkovianRegion$  is partitioned into two non-trivial subintervals $J = \hinterval{i}{\ang} = L_{\ang}^i \cup R_{\ang}^i$, where $L_{\ang}^i := (i,d_{\ang}] \times \{\ang\} $ and $R_{\ang}^i := [d_{\ang},i+1)\times \{\ang\}$, and $d_{\ang}$ is the first coordinate of the intersection point of the image of the curve $\delta$ with the interval $\hinterval{i}{\ang}$. 
%For $|\ang| \geq \pi/6$, one of these two intervals is trivial and the other coincides with $\hinterval{i}{\ang} $, depending on the sign of $\ang$.
%
%
So, the set $\partition = \{ L_\ang^i, \, R_\ang^i : i \in \ZZ_3, |\ang|< \textfrac{\pi}{6}\}$  is a refinement of $\mathcal{J}_{0}$, such that any atom of $\partition$ is \emph{almost} contained in some continuity domain of $T_\lambda$. 
Precisely, given any $S \in \mathcal{J}$, if we denote by $\text{int}(S)$ the interval $S$ without its extreme points, then we have that $\text{int}(S)$ is contained in some continuity region of $T_\lambda$.

The Markovian property of the partition $\MarkovianRegion$ is the following: for any $S\in \partition$, the restriction of $T_\lambda|_{\text{int}(S)}$ extends to a continuous and surjective map $f_{S}:S \to J$, for some $J \in \mathcal{J}$.
%
%Hence, the Markovian property is the following:  for each $S \in \partition$, the restriction $T_\lambda|_{\text{int}(S)}$ is linear, and hence it extends to a continuous and surjective map $f_{S}:S \to J$, for some $J \in \mathcal{J}$. 
In fact, each of these maps is a linear projection between sides of the triangle, along some fixed direction, and hence affine. Furthermore, for any $S \in \partition$ the derivative is constant and  can be computed as 
$(f_{S})' =|S|^{-1}  >1$,
where $|S|$ denotes the length of the interval $S$.  Note that $(f_{J(i,\ang)})'  \to 1$ as $\ang \to \textfrac{\pi}{6}$. These properties do not depend on $\lambda$.

% \vspace{8 pt}
%Consider the maximal forward invariant set $\maximalinvariantset = \cap_{n\geq0}T_{\lambda}(M) \subset M$. 
If we fix some $0 < \lambda < \textfrac{1}{3}$, then the set $\maximalinvariantset$ is properly contained in the Markovian region. That is, there is $\epsilon >0$ such that
$$\maximalinvariantset \subset (0,3) \times (-\textfrac{\pi}{6} + \epsilon, \textfrac{\pi}{6} - \epsilon ) \subset \bigcup_{J \in \mathcal{J}_{0}}J,$$
where equality is excluded.

%In fact, this is true since the largest possible incoming angle is $\pi/2$ the largest outgoing angle is $\lambda \textfrac{\pi}{2} < \textfrac{\pi}{6}$. 
So for any $S \in \partition$ we have that $f_{S}(S) \in \MarkovianRegion$. Moreover, there is a uniform expansion along the horizontal direction: in fact, there is $\rho = \rho(\lambda) >1$ such that $(f_{S})' > \rho$, for any $S \subset \hinterval{i}{\ang}$ such that $|\ang| < \lambda \textfrac{\pi}{2}$.

These properties allow us to 
describe more precisely the dynamics of the pinbilliard map in terms of symbolic dynamics. Let us first state the following lemma:

\begin{lemma} \label{ConjutationPartOne}
Fix some $0< \lambda < \textfrac{1}{3}$. For each $(i, \theta) \in   \ZZ_{3} \times C(\lambda) $ there is a continuous map 
$P(i,\theta): \elshift \to J(i, \theta) \in \MarkovianRegion$
such that $\itin \circ P  = \text{Id}$.
\end{lemma}

%\begin{lemma} %\label{ConjutationPartOne}
%Let $0< \lambda < 1/3$. For any $\theta \in C(\lambda)$, $i \in \ZZ_{3}$, and $w \in \elshift$, there is a point $P \in J(i,\theta) \in \MarkovianRegion$ such that $\itin(P) = w$. Moreover, $P \in \maximalinvariantset$.
%This point is unique, except for those points $w \in \elshift$ whose tail is repeted periodically $(\ldots, -1,1,-1, \ldots )$, which correspond to points that eventually lands on a vertex.
%\end{lemma}

\proof Fix some $0< \lambda < \textfrac{1}{3}$, and consider any$(i, \theta) \in   \ZZ_{3} \times C(\lambda) $. Let $J_{0} = \hinterval{i}{\theta} \in \MarkovianRegion$. 
In order to define the map $P(i,\theta):\elshift \to J_{0}$, consider
a point $w=(w_{n}) \in \elshift$. We take the first $k$ letters of $w$ and  will find a closed interval $H_{k}=H(w_{1}, \ldots, w_{k}) \subset J_{0}$ such that, for any $Q \in H_{k}$, the first $n$ coordinates of $\itin(Q) \in \elshift$ are precisely $w_{1}, \ldots, w_{k}$. 
If $\cap_{k}H_{k} \subset J_{0}$ is a unique point, then the map will then be well-defined.

To define these intervals we shall proceed inductively on $k$:
Recall that $J_0 \in \MarkovianRegion$, and so it splits into two subintervals, $J_{0} = L \cup R$, for some $L,\, R  \in \partition$; then let $H_{1} = R$ if $w_{1}=1$, or $H_{1} = L$ if  $w_{1}=-1$. It is easy to see that the first coordinate of $\itin(Q)$ is $w_{1}$, for any $Q \in H_{1}$. Notice that $f_{J_{0}}(H_{1}) =: J_{1} \in \MarkovianRegion$. 
%Now, for the second letter, since $J_{1} \in \partition$ $J_{1}= L_{1} \cup R_{1}$ for some  $L_{1},\, R_{1}  \in \partition$. Let $H_{2}: = (f_{J})^{-1}(L_{1})$ if $w_{2}=-1$ or  $H_{2}: = (f_{J})^{-1}(R_{1})$ if $w_{2}=1$. Notice that $H_{2} \subset J_{0}$, that  $J_{2}:= f_{J_{1}}(H_{2}) \in \MarkovianRegion$, and that the first two coordinates of $\itin(P)$ are $(w_{1},w_{2})$, for any $P \in H_{2}$. 
Now let us assume that $H_{k-1} \subset J_{0}$ is defined and satisfies 
$$J_{k}  := f_{J_{k-1}}\circ \cdots \circ f_{J_{1}} \circ f_{J_{0}}(H_{k-1})  \in \MarkovianRegion.$$ 
Then $J_{k} = L_{k} \cup R_{k}$ for some $L_{k},\, R_{k}  \in \partition$, and so, in the same way, we can define $H_{k} = (f_{J_{k-1}}\circ \cdots  \circ f_{J_{0}})^{-1}(A)$, where either $A= L_{k}$, or $A= R_{k}$, depending on  $w_{k}$.

We have that $|H_{k}| \to 0$ as $k\to \infty$. In fact, all $J_{k}=\hinterval{\,\cdot\,}{\theta_{k}} \in \MarkovianRegion$ involved in the construction are such that $|\theta_{k}| < \lambda \textfrac{\pi}{2} < \textfrac{\pi}{6}$, and hence $(f_{J_{k}})' > \rho >1$. So $|H_{k}| < \rho ^{-k}$, for any $k \geq 1$, and therefore, there is a unique $P = P(i,\theta)(w) \in \bigcap_{k \geq 1}H_{k}$ and we have $\itin(P) = w$.

This map is continuous by construction and it is injective at all points, except those with a periodic tail (see the definition of the map $\itin$), which correspond to points that eventually hit some vertex of the table, which are mapped to the same points in $J_{0}$.  \qed

As a consequence of the construction in the proof of the previous lemma one obtains the following: on any $J := J(i, \theta) \in \MarkovianRegion$ such that $J \cap \maximalinvariantset \neq \emptyset$, that is $\theta \in C(\lambda)$, any itinerary in $\elshift$ is followed by a unique point in $J$, with exception of those itineraries with periodic tail, which correspond to the same orbit that eventually hits a vertex. In fact, such exceptions correspond to points in the discontinuity set $\DiscontinuityRegion$. Hence, $J \subset \maximalinvariantset$, and this proves that $\maximalinvariantset$ is homeomorphic to $(0,3) \times \elcantor \setminus \DiscontinuityRegion$.

Thus, we can gather together  the maps obtained in Lemma \ref{ConjutationPartOne} to create the following transformation that  depends only on $\lambda$:
$$  P_\lambda: \elshift \times \elshift  \times \ZZ_3 \to \maximalinvariantset;$$
$$ P_\lambda(w, z, i) :=P(i, E_\lambda(z)) (w).$$

Observe that this map is continuous, and if $(p, \ang) \in (0,3) \times C(\lambda)$ then the inverse of $P_\lambda$ is defined by $P_\lambda^{-1}(p,\ang) = (\itin(p,\ang), E_\lambda^{-1}(\ang), \lfloor p \rfloor)$, 
where $\lfloor p \rfloor$ denotes the integer part of $p$.  
Hence, $P_\lambda$ is a homeomorphism in $P_\lambda^{-1}(\maximalinvariantset \setminus \DiscontinuityRegion)$ and conjugates the pinball billiard map to the following symbolic model:
\begin{equation} \label{shiftmodel}
[(w, z),i] \mapsto [(\sigma(w), (w_1, -z_1, -z_2, \ldots )) , i + w_1 ],
\end{equation} 
where $\sigma:\elshift \to \elshift $ is the standard shift, and the addition operation in the last coordinate is in $\ZZ_3$.

This complete description of the dynamics of the pinball map allows us to prove transitivity in $\maximalinvariantset$, as we shall see in the next Lemma.

\begin{lemma}
The pinbilliard map $T_{\lambda}|_{\maximalinvariantset}$ is transitive.
\end{lemma}

\proof Take any two non-empty open subsets $U$ and $V$ of $\maximalinvariantset$. We shall prove that there is some $N \geq 0$ such that $T_{\lambda}^N(U)\cap V \neq \emptyset$. First note that there is some interval $J_{1} := J(i_{1},\theta_{1}) \in \MarkovianRegion$ such that $J_{1} \cap U$ is a non-trivial interval, since $U$ is open. Now, the construction in the proof of Lemma \ref{ConjutationPartOne} allows us to say that there is some interval $H_{1}  \subset J_{1} \cap U$ and $N_{1}\geq 0$ such that $T_{\lambda}^{N_{1}}(H_{1}) = J_{2} \in \MarkovianRegion$. Therefore, there are $(i_{2}, \theta_{2}) \in \ZZ_{3} \times \elcantor$ such that $J_{2}= J(i_{2},  \theta_{2})$. 

On the other hand, there is some interval $J(i_{*}, \tilde \theta ) \in \MarkovianRegion$, such that $J(i_{*}, \tilde \theta ) \cap V \neq \emptyset$ for some $(i_{*}, \tilde \theta) \in \ZZ_{3} \times \elcantor$, 
since $V$ is open.

Now we claim that for any $\epsilon >0$ there is some $N_{2} \geq 0$ and an interval $H_{2} \subset J_{2}$ such that $T_{\lambda}^{N_{2}}(H_{2}) = J_{3} \in \MarkovianRegion$, and hence $J_{3}= J(i_{3},  \theta_{3})$ for some $(i_{3}, \theta_{3}) \in \ZZ_{3} \times \elcantor$, and such that $i_{3} = i_{*}$ and $|\theta_{3} - \tilde \theta| < \epsilon$.
If this claim holds and we take $\epsilon>0$ small enough, it implies that $T_{\lambda}^{N_{1}+N_{2}}(H_{1})\cap V \neq \emptyset$, proving transitivity. 

To get the claim denote $(a_{1}, a_{2},\ldots) = E_{\lambda}^{-1}(\tilde \theta) \in \elshift$.
%$\tilde \theta = \frac{\pi}{3}(a_{1}\lambda +a_{2}\lambda^2+\cdots)$. 
Fix some $k \geq 1$ such that $\lambda^{k+1}< \epsilon$. We will select a point $w \in \elshift$ specifying just the first $k$ coordinates in the following way: let
$$w := ( (-1)^{k-1} a_{k} , (-1)^{k-2}a_{k-1}, \ldots, -a_{2}, a_{1}, w_{k+1}, \ldots) \in \elshift.$$
If we let $z = E_{\lambda}^{-1}(\theta_{2})$, then by the formula in \eqref{shiftmodel} one gets
$$\tau^k(w, z ,i_{2})=(\sigma^k(w) \,, \, (a_{1} , \ldots,a_{k}, z_{1}, \ldots) \,, \, i_{2} + \sum_{s=1}^k (-1)^{s-1} a_{s}).$$

Observe that $| \tilde \theta - E_{\lambda}(a_{1}, \ldots,a_{k}, z_{1}, \ldots) | < \lambda^{k+1} < \epsilon$. So, in order to obtain the claim it is left to prove that $ i_{2} + \sum_{s=1}^k (-1)^{s-1} a_{s} = i_{*} \, \text{mod} \, 3$. However, this is not always true. For instance, let us assume that they differ by 2. Then, to get the claim we may need to change the itinerary $w$, by adjoining two signs $b_{1}$ and $b_{2} \in \{\pm1\}$ at the beginning:
$$w = (b_{2}, b_{1}, (-1)^{k-1} a_{k} , (-1)^{k-2}a_{k-1}, \ldots, -a_{2}, a_{1}, \ldots) \in \elshift,$$
and iterate $k+2$ times to obtain the claim. In fact,
\begin{eqnarray*}
\tau^{k+2}(w, z ,i_{2})=\left( \sigma^{k+2}(w) \,, \, (a_{1} , \ldots,a_{k}, b_{1},b_{2}, z_{1}, \ldots) \right., \\ \left. i_{2} + \sum_{s=1}^k (-1)^{s-1} a_{s} - b_{2} + b_{1} \right).
\end{eqnarray*}

If these numbers differ only by $1$, then the same argument applies by adjoining only one $b_{1} \in \{\pm1\}$. 
\qed

This completes the proof of the topological part of Theorem \ref{MainTheorem}. 

\newcommand{\lamedida}{\mu}

\subsection{Invariant measure on $\limitset$}
The space $ \elshift \times \elshift \times \ZZ_3$ has a natural probability measure $\lamedida$, defined on Borel sets, which is invariant under the map in \eqref{shiftmodel}.  This measure is in fact a Bernoulli measure, and is a product $\lamedida= \nu \times \nu$, where $\nu$ is the invariant Bernoulli measure for the standard shift $\sigma:\elshift \to \elshift$, weighted in order that
$\lamedida( \elshift \times \elshift \times \{i\}) = \textfrac{1}{3}$, for $i \in \ZZ_3$.

On the other hand, one can define a measure $\mu_{\lambda}$ on Borel sets of $\RR$ using the maps $E_{\lambda}$, for $0<\lambda <1$, in the following way: given any measurable set $A \subset \RR$, let 
$$\mu_\lambda(A) := \nu(E_\lambda^{-1}(A)). $$
 
Properties of the measure $\mu_\lambda$, in particular for which values it is
absolutely continuous with respect to Lebesgue measure, have been studied since 1920 by Khintchine and Kolmogoroff \cite{KK25}, and several interesting results have been developed by Erd\"os, Garsia, Pollicott and Solomyak during the 20th century; see \cite{Erd40}, \cite{Sol95}. 
In particular, Jensen and Wintner in \cite{JW35} proved that the measure $\mu_\lambda$ is either purely singular or absolutely continuous with respect to Lebesgue measure in the interval $\frac{\pi}{3} \left[\frac{-\lambda}{1-\lambda}, \frac{\lambda}{1-\lambda} \right]$. If $\lambda < \textfrac{1}{2}$ then the measure $\mu_\lambda$ is purely singular, since it is supported on the Cantor set $\elcantor$.  

So, if $0 < \lambda < \textfrac{1}{3}$, the measure $\lamedida$ gives us an invariant measure for $T_\lambda$ on $\limitset$: Let $m$ be a measure on Borel sets of $A\subset M$ given by
$$m(A) := \lamedida(P_\lambda^{-1}( A \cap \limitset)).$$

\begin{theorem} \label{ergodic}
If $0 < \lambda < \textfrac{1}{3}$, then the measure $m$ is an invariant ergodic measure for the pinball billiard map $T_\lambda$. Moreover, the conditional measures $m|_J$, for $J\in \MarkovianRegion$ such that $J \cap \limitset \neq \emptyset$, are absolutely continuous with respect to the Lebesgue measure on the interval.
\end{theorem}

\proof By construction, the measure $m$ is $T_{\lambda}$-invariant, since it is topologically conjugate to the map in \eqref{shiftmodel}, and hence it inherits all its ergodic properties. On the other hand, conditional measure on horizontal lines are absolutely continuous with respect to the Lebesgue measure, since the map $T_\lambda$ restricted to $J\in \MarkovianRegion$ that $J \cap \limitset \neq \emptyset$ is uniformly expanding with a factor of $\rho>1$. Standard methods \cite[Theorem~2.1]{DMVS}\ complete the proof. \qed

Theorem \ref{ergodic} proves the ergodic part of Theorem \ref{MainTheorem}, and we are done. \qed

\section{Conclusions}

We have built up a detailed picture of the dynamics of non-elastic pinbilliards in an equilateral triangle as the contraction parameter $\lambda$ varies between $0$ and $1$. % $\lambda = 0$ $\lambda = 1$, where the classical elastic billiard is recovered.
Away from these limits, the pinbilliard has dominated splitting, and 
we have established, rigorously for $\lambda < \textfrac{1}{3}$, an intricate limiting dynamics consisting of chaotic strange attractors which evolve as $\lambda$ varies, and which undergo qualitative changes at certain well-defined parameter values.

%These systems exhibit interesting, physically-motivated and calculable properties. 

Preliminary results indicate that non-equilateral triangle tables exhibit similar phenomena, except that the limiting elastic dynamics is then expected to be generically mixing, which seems to be reflected in the attractors found when $\lambda$ close to $1$, which fill out more and more of phase space as $\lambda \to 1$, instead of the shrinking and non-transitive behaviour we have discussed here.
% 
% We have also begun to investigate other regular polygonal billiard tables. Although the dynamics is more complicated, the fundamental phenomena appear to be the same, except that in polygons with parallel edges there are families of attracting period-$2$ orbits which bounce between opposite parallel sides. 
%
Other polygonal billiard tables also seem to exhibit similar behaviour, with additionally the possibility of families of attracting period-$2$ orbits if the table has two parallel sides.

\end{document}